# Efficient forward and inverse uncertainty quantification for dynamical systems based on dimension reduction and Kriging surrogate modeling in functional space

Zhouzhou Song [a, *], Weiyun Xu [b], Marcos A. Valdebenito [a], Matthias G.R. Faes [a, c]

[a] *Chair for Reliability Engineering, TU Dortmund University, 44227 Dortmund, Germany*

[b] *Department of Mechanical Engineering, Tsinghua University, Beijing 100084, China*

[c] *International Joint Research Center for Engineering Reliability and Stochastic Mechanics, Tongji University, Shanghai 200092, China*

\* Corresponding author. E-mail address: zhouzhou.song@tu-dortmund.de

## Abstract

Surrogate models are extensively employed for forward and inverse uncertainty quantification in complex, computation-intensive engineering problems. Nonetheless, constructing high-accuracy surrogate models for complex dynamical systems with limited training samples continues to be a challenge, as capturing the variability in high-dimensional dynamical system responses with a small training set is inherently difficult. This study introduces an efficient Kriging modeling framework based on functional dimension reduction (KFDR) for conducting forward and inverse uncertainty quantification in dynamical systems. By treating the responses of dynamical systems as functions of time, the proposed KFDR method first projects these responses onto a functional space spanned by a set of predefined basis functions, which can deal with noisy data by adding a roughness regularization term. A few key latent functions are then identified by solving the functional eigenequation, mapping the time-variant responses into a low-dimensional latent functional space. Subsequently, Kriging surrogate models with noise terms are constructed in the latent space. With an inverse mapping established from the latent space to the original output space, the proposed approach enables accurate and efficient predictions for dynamical systems. Finally, the surrogate model derived from KFDR is directly utilized for efficient forward and inverse uncertainty quantification of the dynamical system. Through three numerical examples, the proposed method demonstrates its ability to construct highly accurate surrogate models and perform uncertainty quantification for dynamical systems accurately and efficiently.

# 1. Introduction

Dynamical systems are widely encountered in engineering and applied sciences, such as vibratory mechanical systems [1], civil infrastructure [2], and physical or chemical processes [3]. In practice, the performance of a dynamical system is influenced by various uncertainties arising from materials, manufacturing, external forces, and the environment [4-7]. Quantifying the effects of these uncertainties on the system response is crucial. Forward uncertainty quantification and inverse uncertainty quantification are two essential aspects of uncertainty quantification (UQ). Forward UQ focuses on evaluating the uncertainty in system responses caused by uncertain inputs, whereas inverse UQ aims to estimate input uncertainties using observed response data. However, forward UQ is typically conducted using the Monte Carlo sampling method, and inverse UQ often relies on Markov Chain Monte Carlo method within a Bayesian framework [8, 9], both of which require numerous dynamical system simulations. This makes forward UQ and inverse UQ highly inefficient, particularly for computationally expensive problems.

Therefore, surrogate models are widely used in forward UQ and inverse UQ to create computationally efficient models for analysis. Over recent decades, various surrogate modeling approaches have been proposed for emulating dynamical systems. Based on differences in modeling forms, these approaches can be broadly classified into two categories: autoregressive model-based methods and output feature mapping-based methods. Autoregressive model-based methods estimate time-variant responses using past observations or prior predictions. The autoregressive integrated moving average (ARIMA) model [10] is a well-known autoregressive approach that has achieved significant success in time series prediction. However, ARIMA assumes a linear relationship between past history and future forecasts, limiting its applicability to dynamical systems, which often exhibit nonlinearities. Therefore, nonlinear autoregressive models with exogenous input (NARX) [11] were introduced for dynamical systems. These models utilize exogenous inputs, such as time-variant excitation forces, and capture nonlinear relationships between inputs and outputs to achieve higher predictive accuracy. The NARX model enables integration with powerful and widely used surrogate models, such as support vector regression [12, 13], polynomial chaos expansion [14, 15], Kriging (or Gaussian processes) [16, 17], and neural networks [18, 19]. However, determining the optimal time lags for both exogenous and autoregressive inputs is challenging [20], and the NARX model has difficulty handling highly nonlinear dynamical problems [21]. To tackle these challenges, a manifold NARX (mNARX) model [22] was recently introduced, where the input is projected onto a problem-specific manifold that better supports the construction of the NARX model. However, the mNARX model relies on additional physical information.

Output feature mapping-based models aim to map the high-dimensional, time-variant response of a dynamical system into a low-dimensional latent space and construct the surrogate model between the inputs and the latent outputs. The most widely used dimensionality reduction technique for feature mapping is principal component analysis (PCA), also referred to as proper orthogonal decomposition or singular value decomposition in various fields of application. For example, Jacquelin et al. [23] proposed a non-intrusive method that combines PCA with polynomial chaos expansion to model random dynamical systems. Additional studies utilizing PCA to reduce the dimensionality of high-dimensional outputs can be found in [7, 24-27]. However, PCA is a linear mapping method and may not effectively extract features when dealing with highly nonlinear problems. Thus, several methods have been proposed to utilize nonlinear dimensionality reduction techniques for extracting output features. Lee and Carlberg utilized deep convolutional autoencoders to map dynamical systems onto nonlinear manifolds for the purpose of model reduction [28]. Simpson et al. [29] proposed to use autoencoders to infer a latent output space of nonlinear dynamical systems. A systematic investigation of different methods for identifying linear or nonlinear output latent spaces for high-dimensional outputs within the framework of surrogate modeling can be found in [30]. However, accurately identifying the nonlinear latent output space requires a large number of samples, limiting its applicability to problems that involve costly experiments or simulations for generating training samples.

To enhance flexibility and accuracy in inferring the latent output space under noisy conditions and with limited training data, we propose a Kriging modeling framework based on functional dimension reduction (KFDR) for constructing surrogate models for forward and inverse uncertainty quantification in dynamical systems. Fig. 1 presents an overview of the proposed KFDR method. First, instead of viewing the responses of dynamical systems as high-dimensional vectors, we reconsider them from a functional perspective and treat them as functions defined over a specific time interval. From this perspective, we project the time-variant responses onto a functional space spanned by a set of predefined basis functions, which can naturally address noisy data by adding a roughness regularization term. Subsequently, by solving the functional eigenequation, we can capture the majority of variations in the response of the dynamical system through key features in the functional space. The time-variant responses can then be represented as linear combinations of these key latent functions. Thus, the response of the dynamical system is mapped into a low-dimensional latent functional space, with an inverse mapping defined from the latent space to the original output space. Furthermore, Kriging surrogate models with noise terms are constructed in the latent space to account for errors arising from limited data and feature mapping, enabling accurate and efficient predictions of dynamical systems. Finally, the surrogate model constructed using KFDR is directly employed for efficient forward and Bayesian inverse UQ of the dynamical system.

The remainder of this paper is organized as follows. Section 2 introduces the fundamentals of forward and inverse UQ approaches for dynamical systems. Section 3 outlines the details of the proposed KFDR method. Section 4 presents case studies and discusses their results. Finally, Section 5 concludes the paper and suggests potential future research directions.

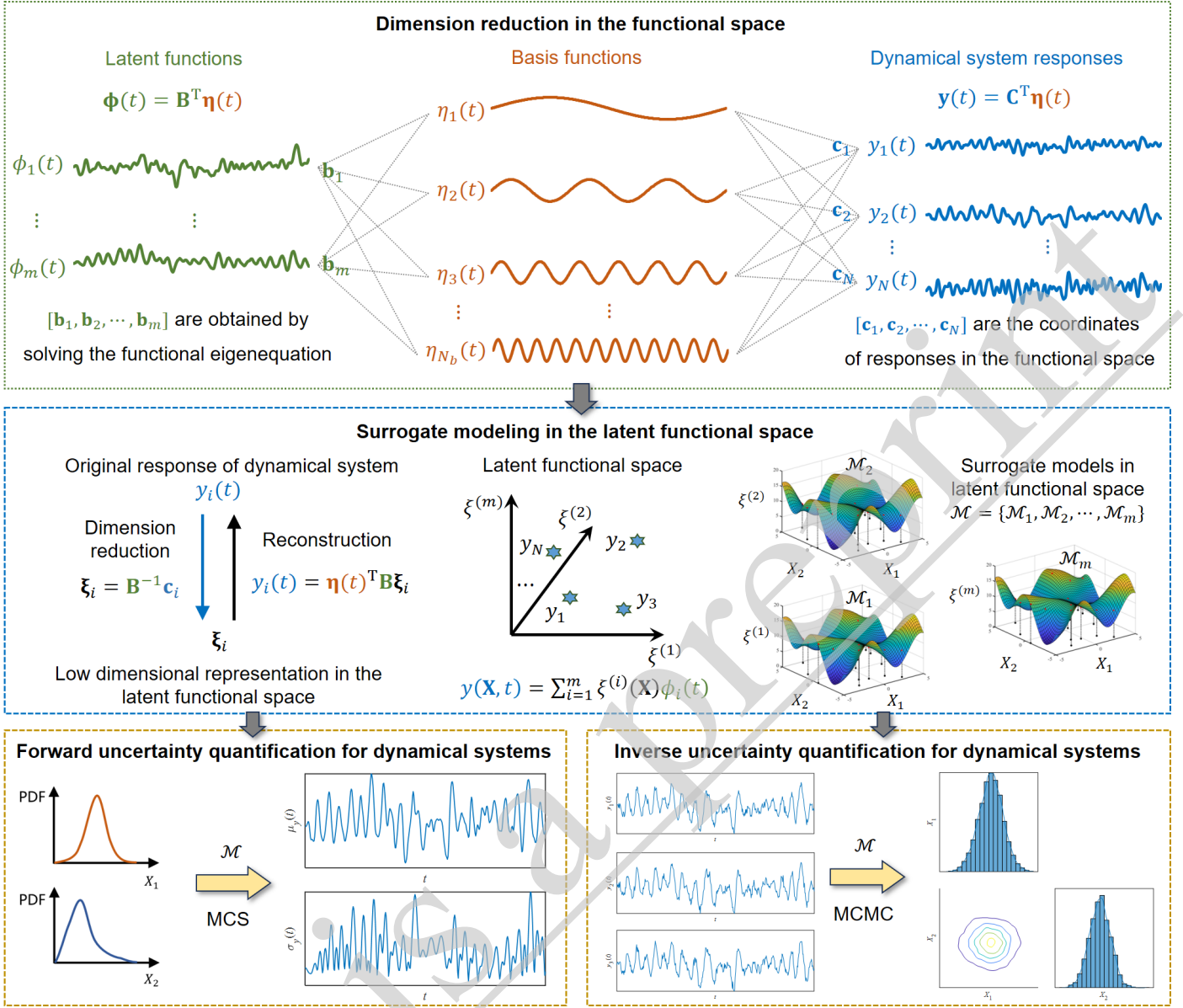

**Fig. 1**. Overview of the proposed framework.

## 2. Problem statement

A response of interest of a dynamical system can be expressed as $Y(\mathbf{X}, t), t \in [t_0, t_e]$, where $\mathbf{X} = [X_1, \cdots, X_p]^\mathrm{T} \in \mathbb{R}^p$ is the input vector. The purpose of forward uncertainty quantification is to obtain statistical information about the time-variant output $Y$ given the probability density function of the input $\mathbf{X} \sim f_\mathbf{X}(\mathbf{x})$. The statistics of interest typically include the mean function,

$$\mu_Y(t) = \int Y(\mathbf{X}, t) f_\mathbf{X}(\mathbf{x}) \mathrm{d}\mathbf{x}, \tag{1}$$

the standard deviation function,

$$\sigma_Y(t) = \sqrt{\int (Y(\mathbf{X}, t) - \mu_Y(t))^2 f_\mathbf{X}(\mathbf{x}) \mathrm{d}\mathbf{x}}, \tag{2}$$

and the probability density function of the output at different time nodes. For reliability analysis, the distributions of the maximum or minimum values over a specified time interval are also of interest. Since $Y(\mathbf{X}, t)$ often lacks an analytic expression, it is generally intractable to compute Eq. (1) and Eq. (2) directly. As a result, simulation methods are commonly used for forward UQ, with Monte Carlo simulation being one of the most widely used approaches. However, these methods require numerous evaluations of $Y(\mathbf{X}, t)$ to obtain precise results, which is computationally prohibitive especially for engineering applications that rely on costly simulations. To reduce the computational burden, a surrogate model of $Y(\mathbf{X}, t)$ needs to be constructed.

Forward UQ relies on the input uncertainty information $f_\mathbf{X}(\mathbf{x})$ to obtain the uncertainty information of outputs. However, obtaining accurate $f_\mathbf{X}(\mathbf{x})$ is often challenging in engineering applications, as it may require a large number of experiments. In some cases, prior knowledge can be used to determine $f_\mathbf{X}(\mathbf{x})$, but this approach can be subjective and may lead to inaccurate forward UQ results. In this case, inverse UQ is needed to infer the uncertainty of the input based on observed response data. Inverse UQ is typically based on a Bayesian framework [31]. First, the input parameters are assumed to follow certain prior distributions, which are then updated according to the observed response data to obtain the posterior distributions, ensuring that the simulation results are consistent with the response data. Markov Chain Monte Carlo sampling is commonly used to compute the posterior distributions, which requires numerous evaluations of $Y(\mathbf{X}, t)$. Therefore, a computationally efficient surrogate model is needed for effective inverse UQ.

Forward UQ and inverse UQ are two essential components of uncertainty quantification for dynamical systems. However, both forward UQ and inverse UQ require numerous system evaluations to obtain responses, making them computationally inefficient for complex problems. To address this, this paper proposes an efficient surrogate-based forward UQ and inverse UQ framework for dynamical systems.

# 3. Methodology

As described in Section 2, the key to efficient forward UQ and inverse UQ of a dynamical system is constructing a surrogate model of it. In this section, we first introduce how to represent the responses of dynamical systems from a functional perspective. Then, dimension reduction and surrogate modeling are performed in the functional space. Subsequently, surrogate-based forward and inverse uncertainty quantification are described at the end of this section.

## 3.1. Dimension reduction in functional space

Since the output of a dynamical system is a function of time, treating it from a functional perspective allow us to obtain more useful information than traditional linear dimensionality reduction methods. For a square-integrable stochastic process $Y(t), t \in [t_0, t_e]$, let $\mu(t) = \mathbb{E}(Y(t))$ be the mean function of $Y$, and let $Y^c(t) = Y(t) - \mu(t)$ be the centered stochastic process. The covariance function of $Y$ is defined as:

$$c(s,t) = \text{Cov}(Y(s), Y(t)) = \mathbb{E}[Y^c(s)Y^c(t)]. \tag{3}$$

Given that the covariance function is symmetric and positive semi-definite, Mercer's theorem [32] implies that:

$$c(s,t) = \sum_{k=1}^{\infty} \lambda_k \phi_k(s) \phi_k(t), \tag{4}$$

where $\lambda_1 \geq \lambda_2 \geq \cdots \geq 0$ are the eigenvalues and $\phi_1, \phi_2, \cdots$ are the corresponding orthonormal eigenfunctions of the covariance operator:

$$C: L^2([t_0, t_e]) \to L^2([t_0, t_e]), C[f](t) = \int_{t_0}^{t_e} c(s,t) f(s) ds, \tag{5}$$

where $L^2([t_0, t_e])$ refers to the space of square-integrable functions defined on $[t_0, t_e]$. Then, by Karhunen-Loève expansion, we have:

$$Y(t) = \mu(t) + \sum_{k=1}^{\infty} \xi^{(k)} \phi_k(t), \tag{6}$$

where

$$\xi^{(k)} = \langle Y^c, \phi_k \rangle = \int_{t_0}^{t_e} Y^c(t) \phi_k(t) dt, k = 1, 2, \ldots \tag{7}$$

are uncorrelated random variables with zero mean and variances of $\lambda_1, \lambda_2, \cdots$, respectively. $\xi^{(k)}$ is the principal component score associated with the *k*-th eigenfunction $\phi_k$ and is the projection of $Y^c(t)$ in the direction of the *k*-th eigenfunction $\phi_k$.

The eigen functions $\phi_1, \phi_2, \cdots$ can be obtained by solving the Fredholm integral equation of the second kind, expressed as:

$$\int_{t_0}^{t_e} c(s,t)\phi(t)dt = \lambda\phi(s). \tag{8}$$

In practice, the continuous eigenproblem in Eq. (8) is discretized into a matrix eigenproblem to facilitate the solution of the integral equation. This is achieved by projecting the dynamical system response $Y(t)$ and the eigen function $\phi(t)$ onto a functional space spanned by predefined basis functions. The covariance function $c(s,t)$ is then estimated using samples of $Y(t)$. Given a training data set with $N$ samples, $\mathcal{D} = \{(\mathbf{x}_i, \mathbf{y}_i), i = 1,2,\cdots,N\}$, where $\mathbf{y}_i$ is an $N_t \times 1$ output vector, and $N_t$ is the number of discretized time nodes. Each $\mathbf{y}_i$ represents a response function $Y_i(t)$. First, the output data is centered as:

$$\mathbf{y}_i^c = \mathbf{y}_i - \frac{1}{N}\sum_{k=1}^{N} \mathbf{y}_k. \tag{9}$$

Then, the centered time-variant output functions $\{Y_1^c(t), Y_2^c(t), \cdots, Y_N^c(t)\}$ are expressed as linear combinations of predefined basis functions $\{\eta_1(t), \eta_2(t), \cdots, \eta_{N_b}(t)\}$ as:

$$\mathbf{Y}^c(t) = [Y_1^c(t), Y_2^c(t), \cdots, Y_N^c(t)]^\mathrm{T} = \mathbf{C}^\mathrm{T}\boldsymbol{\eta}(t). \tag{10}$$

where $N_b$ is the number of basis functions, $\mathbf{C} = [\mathbf{c}_1, \mathbf{c}_2, \cdots, \mathbf{c}_N]$ is the $N_b \times N$ coefficient matrix, $\mathbf{c}_i, i = 1,2,\cdots,N$ are $N_b \times 1$ coefficient vectors, and $\boldsymbol{\eta}(t) = [\eta_1(t), \eta_2(t), \cdots, \eta_{N_b}(t)]^\mathrm{T}$. $\mathbf{c}_i$ can be obtained by minimizing the sum of squared error between the observed and estimated response:

$$\mathbf{c}_i = \underset{\mathbf{c}}{\mathrm{argmin}} \sum_{j=1}^{N_t} \left[y_{ij}^c - Y_i^c(t_j)\right]^2 = \underset{\mathbf{c}}{\mathrm{argmin}}(\mathbf{y}_i^c - \mathbf{Hc})^\mathrm{T}(\mathbf{y}_i^c - \mathbf{Hc}), \tag{11}$$

where

$$\mathbf{H} = \begin{bmatrix} \eta_1(t_1) & \eta_2(t_1) & \cdots & \eta_{N_b}(t_1) \\ \eta_1(t_2) & \eta_2(t_2) & \cdots & \eta_{N_b}(t_2) \\ \vdots & \vdots & \ddots & \vdots \\ \eta_1(t_{N_t}) & \eta_2(t_{N_t}) & \cdots & \eta_{N_b}(t_{N_t}) \end{bmatrix} \tag{12}$$

is an $N_t \times N_b$ matrix whose elements correspond to the values of various basis functions at different time nodes.

If the dynamical system responses include noise, such as that arising from measurements, a roughness regularization term [33] can be added to Eq. (11) as follows:

$$\mathbf{c}_i = \underset{\mathbf{c}}{\mathrm{argmin}} \left\{ (\mathbf{y}_i^c - \mathbf{Hc})^\mathrm{T}(\mathbf{y}_i^c - \mathbf{Hc}) + \tau \int_{t_0}^{t_e} [D^2 Y_i^c(t)]^2 dt \right\}, \tag{13}$$

where $D^2 Y_i^c(t)$ is the second derivative of $Y_i^c(t)$, and the integrated squared second derivative measures the roughness of $Y_i^c(t)$. $\tau$ is the smoothing parameter and is non-negative. A large $\tau$ will cause $Y_i^c(t)$ to exhibit minimal fluctuations. As $\tau$ reaches zero, $Y_i^c(t)$ will attempt to pass through each sample point as closely as possible, potentially leading to erratic behavior in certain regions. By substituting $Y_i^c(t) =$

$\boldsymbol{\eta}(t)^{\mathrm{T}}\mathbf{c}_i$ into the roughness penalty term in Eq. (13), we obtain:

$$\tau \int_{t_0}^{t_e}[D^2 Y_i^c(t)]^2 dt = \tau \mathbf{c}_i^T \left[\int_{t_0}^{t_e} D^2\boldsymbol{\eta}(t)D^2\boldsymbol{\eta}(t)^{\mathrm{T}} dt\right]\mathbf{c}_i. \tag{14}$$

Let $\mathbf{R} = \int_{t_0}^{t_e} D^2\boldsymbol{\eta}(t)D^2\boldsymbol{\eta}(t)^{\mathrm{T}} dt$, where $\mathbf{R}$ is an $N_b \times N_b$ symmetric matrix with elements $\mathbf{R}_{ij} = \int_{t_0}^{t_e} D^2\eta_i(t)D^2\eta_j(t) dt$, $i,j = 1,2,\cdots,N_b$. The analytical solution to Eq. (13) is:

$$\mathbf{c}_i = (\mathbf{H}^{\mathrm{T}}\mathbf{H} + \tau\mathbf{R})^{-1}\mathbf{H}^{\mathrm{T}}\mathbf{y}_i^c. \tag{15}$$

Typically, the smoothing parameter $\tau$ can be determined through cross-validation. However, cross-validation is usually computationally expensive. In this research, we employ the generalized cross-validation (GCV) measure [34], which serves as a more efficient alternative to the standard cross-validation procedure. The GCV measure is expressed as:

$$\mathrm{GCV}(\tau) = \frac{N}{[N - \mathrm{trace}(\mathbf{S}(\tau))]^2}\sum_{i=1}^{N}(\mathbf{y}_i^c - \mathbf{H}\mathbf{c}_i)^{\mathrm{T}}(\mathbf{y}_i^c - \mathbf{H}\mathbf{c}_i), \tag{16}$$

where $\mathbf{S}(\tau) = \mathbf{H}(\mathbf{H}^{\mathrm{T}}\mathbf{H} + \tau\mathbf{R})^{-1}\mathbf{H}^{\mathrm{T}}$. Then, the value of $\tau$ that minimizes $\mathrm{GCV}(\tau)$ is selected for use in Eq. (15). In practice, it is not necessary to compute the exact minimum of $\mathrm{GCV}(\tau)$. Instead, a grid search on a logarithmic scale can be performed to find the optimal $\tau$. For example, the range of $\log_{10}\tau$ can be set to [-6, 6] and divided into uniform grids. The GCV value is then calculated for each grid point, and the $\tau$ corresponding to the grid point with the minimum GCV value is selected for use in Eq. (15). Table 1 presents the pseudocode for determining $\tau$.

**Table 1**
Pseudocode of determining the smoothing parameter $\tau$.

| Algorithm 1: Determination of the smoothing parameter $\tau$ |
|---|
| **Input**: centralized time-variant output samples $\{\mathbf{y}_1^c, \mathbf{y}_2^c, \cdots, \mathbf{y}_N^c\}$ and basis functions $\{\eta_1(t), \eta_2(t), \cdots, \eta_{N_b}(t)\}$ |
| **Output**: the smoothing parameter $\tau$ |
| 1:   Generate $N_\tau$ values for $\tau$: $\tau_i \leftarrow 10^{-6+12(i-1)/(N_\tau-1)}, i = 1,2,\cdots,N_\tau$ |
| 2:   Calculate the values of $\mathrm{GCV}(\tau_i)$ for each $\tau_i$ |
| 3:   $\tau \leftarrow \min_{\tau_i \in \{\tau_1,\cdots,\tau_{N_\tau}\}} \mathrm{GCV}(\tau_i)$ |

The commonly used basis functions, $\boldsymbol{\eta}(t)$, include Fourier basis functions and spline basis functions. Among the spline basis functions, B-spline functions [35] are extensively used, particularly in applications such as computer-aided design and computer graphics. Fourier series are well-suited for periodic data and functions, whereas B-spline series are more appropriate for non-periodic data. Fig. 2 illustrates the Fourier and B-spline basis functions over the time interval [0, 1]. Although Fourier basis functions cannot directly represent non-periodic data, the data can be mirrored along the time axis to generate periodic data, as illustrated in Fig. 3. Once periodicity is established, Fourier functions can then be utilized as basis functions.

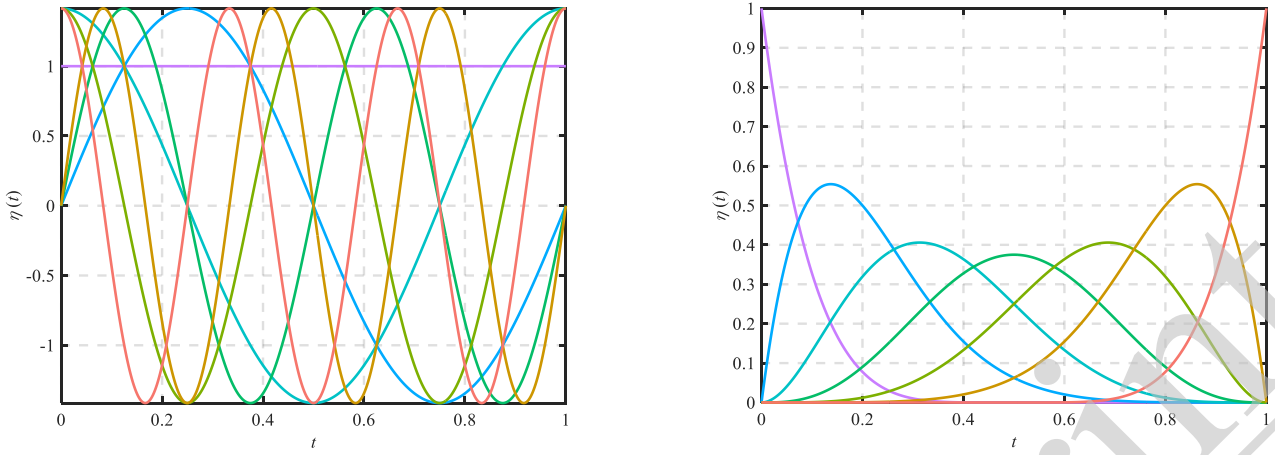

**Fig. 2**. Illustrations of Fourier basis functions (left) and B-spline basis functions (right) over the time interval [0, 1].

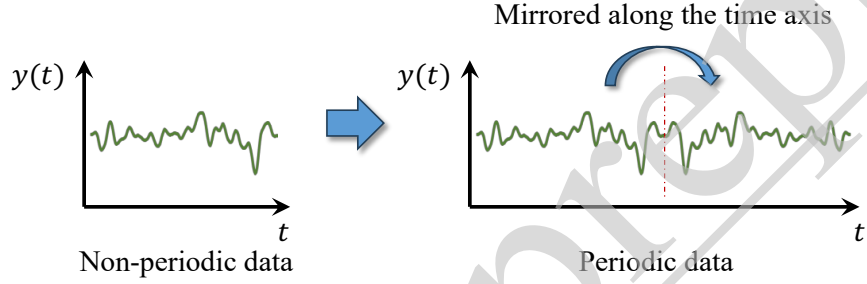

**Fig. 3**. Generate periodic data from non-periodic data.

Determining the number of basis functions $N_b$ is crucial, as it directly influences the representation accuracy. In this study, we develop an error-based approach to select the appropriate $N_b$. After obtaining the coordinates $\mathbf{c}_i$ of $\mathbf{y}_i$ in the functional space using Eq. (15), the error between $\mathbf{y}_i^c$ and $\mathbf{H}\mathbf{c}_i$ indicates the accuracy of the projection from the original time-variant response space to the functional space spanned by $\{\eta_1(t), \eta_2(t), \cdots, \eta_{N_b}(t)\}$. To quantify the deviation between $\mathbf{y}_i^c$ and $\mathbf{H}\mathbf{c}_i$, we use the normalized root mean square error (NRMSE):

$$\mathrm{NRMSE}(\mathbf{y}_i^c, \mathbf{H}\mathbf{c}_i) = \frac{\|\mathbf{y}_i^c - \mathbf{H}\mathbf{c}_i\|_2}{\max \mathbf{y}_i^c - \min \mathbf{y}_i^c}. \tag{17}$$

The average NRMSE of the training set is utilized to quantify the overall error:

$$\delta = \frac{1}{N}\sum_{i=1}^{N} \mathrm{NRMSE}(\mathbf{y}_i^c, \mathbf{H}\mathbf{c}_i). \tag{18}$$

Note that, given a training set and a basis series, $\delta$ depends solely on the number of basis functions, $N_b$. Consequently, $N_b$ can be incrementally increased from a starting value until the relative error between two consecutive $\delta$ values falls below a specified threshold $\delta_r$. $\delta_r$ is set to be 0.05 in this research. The focus on the difference between two consecutive $\delta$ values, rather than $\delta$ itself, arises from the fact that when noise is present, increasing $N_b$ causes $\mathrm{NRMSE}(\mathbf{y}_i^c, \mathbf{H}\mathbf{c}_i)$ to approach a value greater than zero instead of zero. In such cases, using $\delta$ as the convergence criterion may result in failure to converge. Table 2 provides the pseudocode for the error-based approach to determine $N_b$.

**Table 2**
Pseudocode for the error-based approach to determine the number of basis functions $N_b$.

| |
|---|
| **Algorithm 2**: Error-based approach to determine the number of basis functions $N_b$ |
| **Input**: centralized time-variant output samples $\{\mathbf{y}_1^c, \mathbf{y}_2^c, \cdots, \mathbf{y}_N^c\}$ and a basis function system $\{\eta_1(t), \eta_2(t), \cdots\}$ |
| **Output**: number of basis functions $N_b$ and the corresponding $\tau$ |
| 1:   $N_b \leftarrow N_b^0$, where $N_b^0$ is a given positive integer |
| 2:   Determine $\tau$ by using Algorithm 1 |
| 3:   Compute the matrix $\mathbf{H}$: $\mathbf{H}_{ij} \leftarrow \eta_j(t_i), i = 1,2,\cdots,N_t, j = 1,2,\cdots,N_b$ |
| 4:   $\mathbf{c}_i \leftarrow (\mathbf{H}^T\mathbf{H} + \tau\mathbf{R})^{-1}\mathbf{H}^T\mathbf{y}_i^c, i = 1,\cdots,N$ by using the basis functions $\{\eta_1(t), \eta_2(t), \cdots, \eta_{N_b}(t)\}$ |
| 5:   $\delta_1 \leftarrow N^{-1}\sum_{i=1}^{N} \text{NRMSE}(\mathbf{y}_i^c, \mathbf{H}\mathbf{c}_i)$ |
| 6:   $k \leftarrow 1$ |
| 7:   **While** 1 |
| 8:       $N_b \leftarrow N_b + kN_b^0$ |
| 9:       Determine $\tau$ by using Algorithm 1 |
| 10:      Compute the matrix $\mathbf{H}$: $\mathbf{H}_{ij} \leftarrow \eta_j(t_i), i = 1,2,\cdots,N_t, j = 1,2,\cdots,N_b$ |
| 11:      $\mathbf{c}_i \leftarrow (\mathbf{H}^T\mathbf{H} + \tau\mathbf{R})^{-1}\mathbf{H}^T\mathbf{y}_i^c, i = 1,\cdots,N$ by using the basis functions $\{\eta_1(t), \eta_2(t), \cdots, \eta_{N_b}(t)\}$ |
| 12:      $\delta_2 \leftarrow N^{-1}\sum_{i=1}^{N} \text{NRMSE}(\mathbf{y}_i^c, \mathbf{H}\mathbf{c}_i)$ |
| 13:      **If** $|\delta_1 - \delta_2|/\delta_2 < \delta_r$ |
| 14:          **Break** |
| 15:      **End If** |
| 16:      $\delta_1 \leftarrow \delta_2$ |
| 17:      $k \leftarrow k + 1$ |
| 18:   **End While** |

After obtaining the coefficient matrix $\mathbf{C}$, the covariance function $c(s,t)$ is estimated as:

$$c(s,t) = \frac{1}{N-1}\boldsymbol{\eta}(s)^T\mathbf{C}\mathbf{C}^T\boldsymbol{\eta}(t). \tag{19}$$

The eigenfunction $\phi(s)$ in Eq. (8) can also be approximated by the basis functions $\boldsymbol{\eta}(s)$ as:

$$\phi(s) = \mathbf{b}^T\boldsymbol{\eta}(s) = \boldsymbol{\eta}(s)^T\mathbf{b}, \tag{20}$$

where $\mathbf{b}$ is an $N_b \times 1$ vector and stands for the coordinates of $\phi(s)$ in the functional space spanned by $\boldsymbol{\eta}(s)$. Substituting Eq. (19) and Eq. (20) into Eq. (8), we obtain:

$$\int_{t_0}^{t_e} c(s,t)\phi(t)dt = \lambda\boldsymbol{\eta}(s)^T\mathbf{b},$$

$$\int_{t_0}^{t_e} \frac{1}{N-1}\boldsymbol{\eta}(s)^T\mathbf{C}\mathbf{C}^T\boldsymbol{\eta}(t)\boldsymbol{\eta}(t)^T\mathbf{b}dt = \lambda\boldsymbol{\eta}(s)^T\mathbf{b},$$

$$\frac{1}{N-1}\boldsymbol{\eta}(s)^T\mathbf{C}\mathbf{C}^T\left[\int_{t_0}^{t_e}\boldsymbol{\eta}(t)\boldsymbol{\eta}(t)^Tdt\right]\mathbf{b} = \lambda\boldsymbol{\eta}(s)^T\mathbf{b}. \tag{21}$$

Let $\mathbf{W} = \int_{t_0}^{t_e}\boldsymbol{\eta}(t)\boldsymbol{\eta}(t)^Tdt$, where $\mathbf{W}$ is an $N_b \times N_b$ symmetric matrix with elements $\mathbf{W}_{ij} =$ $\int_{t_0}^{t_e} \eta_i(t)\eta_j(t)dt$, $i,j = 1,2,\cdots,N_b$. Then the discrete form of eigenequation is obtained:

$$\frac{1}{N-1}\boldsymbol{\eta}(s)^{\mathrm{T}}\mathbf{CC}^{\mathrm{T}}\mathbf{Wb} = \lambda\boldsymbol{\eta}(s)^{\mathrm{T}}\mathbf{b}. \tag{22}$$

Since this equation must hold for all $s$, we obtain:

$$\frac{1}{N-1}\mathbf{CC}^{\mathrm{T}}\mathbf{Wb} = \lambda\mathbf{b}. \tag{23}$$

By defining $\mathbf{u} = \mathbf{W}^{1/2}\mathbf{b}$, we need to solve finally a symmetric eigenvalue problem:

$$\frac{1}{N-1}\mathbf{W}^{1/2}\mathbf{CC}^{\mathrm{T}}\mathbf{W}^{1/2}\mathbf{u} = \lambda\mathbf{u}, \tag{24}$$

and compute $\mathbf{b} = \mathbf{W}^{-1/2}\mathbf{u}$ for each eigenvector. Note that Fourier basis functions are orthogonal to each other. Consequently, the matrix $\mathbf{W}$ reduces to an identity matrix for Fourier basis functions, and Eq. (23) simplifies to performing standard principal component analysis on the coefficient matrix $\mathbf{C}$. In contrast, since B-spline basis functions are generally not orthogonal, it is necessary to compute $\mathbf{W}$ and solve the eigenproblem in Eq. (24).

In practice, only the first few eigenfunctions $\{\phi_1(t), \phi_2(t), \cdots, \phi_m(t)\}$ are sufficient to represent $Y(t)$. There are several methods to determine the value of $m$, including the variance proportion-based approach [36-38], the Bayesian information criterion-based approach [39], the reconstruction error-based approach [26], and the ladle estimator-based approach [40, 41]. In this study, we adopt the 99% variance proportion-based approach due to its simplicity and efficiency. Specifically, $m$ is chosen as the smallest value that satisfies:

$$\frac{\sum_{i=1}^{m}\lambda_i}{\sum_{i=1}^{N_b}\lambda_i} \geq 99\%, \tag{25}$$

where $\lambda_1, \lambda_2, \cdots, \lambda_{N_b}$ are the eigenvalues of $(N-1)^{-1}\mathbf{W}^{1/2}\mathbf{CC}^{\mathrm{T}}\mathbf{W}^{1/2}$.

Once the eigenfunctions are obtained, the original high-dimensional time-variant response can be reduced to a low-dimensional vector, and the original response can be reconstructed from its low-dimensional representation. Let $\mathbf{B} = [\mathbf{b}_1, \mathbf{b}_2, \cdots, \mathbf{b}_m]$, which is an $N_b \times m$ matrix. For a new time-variant response $\mathbf{y}^*$, its low-dimensional representation $\boldsymbol{\xi}^*$ can be obtained in the same way as in Eq. (13) and Eq. (15):

$$\boldsymbol{\xi}^* = \underset{\boldsymbol{\xi}}{\mathrm{argmin}}\{(\mathbf{y}^* - \mathbf{HB}\boldsymbol{\xi})^{\mathrm{T}}(\mathbf{y}^* - \mathbf{HB}\boldsymbol{\xi}) + \tau\boldsymbol{\xi}^{\mathrm{T}}\mathbf{B}^{\mathrm{T}}\mathbf{RB}\boldsymbol{\xi}\} = \mathbf{B}^{-1}(\mathbf{H}^{\mathrm{T}}\mathbf{H} + \tau\mathbf{R})^{-1}\mathbf{H}^{\mathrm{T}}\mathbf{y}^*. \tag{26}$$

For the training samples $\{\mathbf{y}_1, \mathbf{y}_2, \cdots, \mathbf{y}_N\}$, their low-dimensional representation can be directly obtained as $\boldsymbol{\xi}_i = \mathbf{B}^{-1}\mathbf{c}_i$, for $i = 1, \cdots, N$, where $\mathbf{c}_i$ is computed as in Eq. (15) and $\boldsymbol{\xi}_i$ is an $m \times 1$ vector. Additionally, for a low-dimensional vector $\hat{\boldsymbol{\xi}}$ in the latent functional space, the time-variant response is reconstructed as $\hat{y}(t) = \boldsymbol{\eta}(t)^{\mathrm{T}}\mathbf{B}\hat{\boldsymbol{\xi}}$. By performing dimension reduction in the functional space, we connect the high-dimensional time-variant response to a low-dimensional vector in the latent functional space. Therefore, we can construct surrogate models between the inputs and the latent outputs to predict the original time-variant response. Table 3 presents the pseudocode for performing dimension reduction in the functional space.

**Table 3**
Pseudocode of dimension reduction in the functional space for the time-variant response.

| Algorithm 3: Dimension reduction in the functional space |
|---|
| **Input**: time-variant output samples $\{\mathbf{y}_1, \mathbf{y}_2, \cdots, \mathbf{y}_N\}$ |
| **Output**: low-dimensional representations for the output samples $\{\boldsymbol{\xi}_1, \boldsymbol{\xi}_2, \cdots, \boldsymbol{\xi}_N\}$ |
| 1: Centralize the data $\mathbf{y}_i^c \leftarrow \mathbf{y}_i - N^{-1}\sum_{i=1}^{N} \mathbf{y}_i, i = 1,2,\cdots,N$ |
| 2: Select a basis function system $\{\eta_1(t), \eta_2(t), \cdots\}$ and determine $N_b$ and $\tau$ by Algorithm 2 |
| 3: Compute the matrix $\mathbf{H}$: $\mathbf{H}_{ij} \leftarrow \eta_j(t_i), i = 1,2,\cdots,N_t, j = 1,2,\cdots,N_b$ |
| 4: $\mathbf{c}_i \leftarrow (\mathbf{H}^T\mathbf{H} + \tau\mathbf{R})^{-1}\mathbf{H}^T\mathbf{y}_i^c, i = 1,\cdots,N$ and $\mathbf{C} \leftarrow [\mathbf{c}_1, \mathbf{c}_2, \cdots, \mathbf{c}_N]$ |
| 5: $\mathbf{W} \leftarrow \int_{t_0}^{t_e} \boldsymbol{\eta}(t)\boldsymbol{\eta}(t)^T dt$ |
| 6: Solve the symmetric eigenvalue problem $(N-1)^{-1}\mathbf{W}^{1/2}\mathbf{C}\mathbf{C}^T\mathbf{W}^{1/2}\mathbf{u} = \lambda\mathbf{u}$ |
| 7: Determine $m$ with the 99% variance proportion criterion and obtain the retained eigen pairs: $\{\lambda_1, \mathbf{u}_1\}, \cdots, \{\lambda_m, \mathbf{u}_m\}$ |
| 8: $\mathbf{b}_k \leftarrow \mathbf{W}^{-1/2}\mathbf{u}_k, k = 1,2,\cdots,m$ and $\mathbf{B} \leftarrow [\mathbf{b}_1, \mathbf{b}_2, \cdots, \mathbf{b}_m]$ |
| 9: $\boldsymbol{\xi}_i \leftarrow \mathbf{B}^{-1}\mathbf{c}_i, i = 1,\cdots,N$ |

### 3.2. Kriging-based emulator for learning dynamical systems

After performing dimension reduction in the functional space, time-variant response $\mathbf{y}$ is reduced to an $m \times 1$ vector $\boldsymbol{\xi} = [\xi^{(1)}, \xi^{(2)}, \cdots, \xi^{(m)}]^T$ in the latent functional space. Since the latent functions $\phi_k(t)$ are orthogonal to each other, all $\xi^{(k)}, k = 1,2,\cdots,m$ are uncorrelated. Therefore, to emulate the dynamical system, we can construct a surrogate model for each $\xi^{(j)}$ with respect to the input $\mathbf{X}$ and use these models to predict the system's response. In this study, we use the Kriging surrogate modeling method due to its ability to quantify model prediction uncertainty, a highly valuable feature for assessing the surrogate model's quality or supporting active learning. The training data $\mathcal{D} = \{(\mathbf{x}_i, \boldsymbol{\xi}_i), i = 1,2,\cdots,N\}$ in the latent space may contain noise due to limited data and dimensionality reduction. Therefore, the Ordinary Kriging model with noise term is used for surrogate modeling in the latent space:

$$\xi(\boldsymbol{x}) = \mu + Z(\boldsymbol{x}) + \varepsilon, \tag{27}$$

where $\mu$ is the global mean, $Z(\mathbf{x}) \sim \mathrm{GP}(0, k(\mathbf{x}, \mathbf{x}'))$ is a zero mean Gaussian process, $\varepsilon$ is a zero-mean Gaussian noise with covariance matrix $\boldsymbol{\Sigma}_n$. This paper assumes homoscedastic noise, where $\boldsymbol{\Sigma}_n = \sigma_n^2 \mathbf{I}$ (**I** being the identity matrix). It is worth noting that the assumption of homoscedastic noise can be relaxed to accommodate heteroscedastic noise, but such considerations are beyond the scope of this work. $k(\mathbf{x}, \mathbf{x}') = \mathbb{E}[Z(\mathbf{x})Z(\mathbf{x}')]$ is the covariance function (or kernel function) of $Z(\mathbf{x})$. Among numerous existing kernel functions, the Gaussian kernel function is commonly used:

$$k(\mathbf{x}, \mathbf{x}') = \sigma_Z^2 \exp\{-(\mathbf{x} - \mathbf{x}')^T \boldsymbol{\Theta}(\mathbf{x} - \mathbf{x}')\}, \tag{28}$$

where $\sigma_Z^2$ is the variance of $Z(\mathbf{x})$, $\boldsymbol{\Theta} = \mathrm{diag}(\boldsymbol{\theta})$ and $\boldsymbol{\theta} = [\theta_1, \theta_2, \cdots \theta_p]^T$ are scaling parameters to characterize the variability of the Gaussian process.

Given a training data set $\mathcal{D}_j = \{(\mathbf{x}_i, \xi_i^{(j)}), i = 1,2,\cdots,N\}$ for the $j$-th component of $\boldsymbol{\xi}$, $\mu$, $\sigma_Z^2$, $\boldsymbol{\theta}$ and $\sigma_n^2$ are obtained by maximizing the marginal log likelihood function as follows:

$$\hat{\mu}, \hat{\sigma}_Z^2, \widehat{\boldsymbol{\theta}}, \hat{\sigma}_n^2 = \underset{\mu, \sigma_Z^2, \boldsymbol{\theta}, \sigma_n^2}{\operatorname{argmax}} \ln(\xi|\mathbf{x}, \mu, \sigma_Z^2, \boldsymbol{\theta}, \sigma_n^2),$$

$$\ln(\xi|\mathbf{x}, \mu, \sigma_Z^2, \boldsymbol{\theta}, \sigma_n^2) = -\frac{1}{2}\big(\boldsymbol{\xi}^{(j)} - \mathbf{1}\mu\big)^{\mathrm{T}} (K_{\mathbf{xx}} + \sigma_n^2 \mathbf{I}_N)^{-1} \big(\boldsymbol{\xi}^{(j)} - \mathbf{1}\mu\big) - \frac{1}{2}\ln|K_{\mathbf{xx}} + \sigma_n^2 \mathbf{I}_N| - \frac{N}{2}\ln 2\pi, \quad (29)$$

where $\boldsymbol{\xi}^{(j)} = [\xi_1^{(j)}, \xi_2^{(j)}, \cdots, \xi_N^{(j)}]^{\mathrm{T}}$, $\mathbf{1}$ is an $N \times 1$ vector of ones, $K_{\mathbf{xx}}$ is the $N \times N$ covariance matrix with $(K_{\mathbf{xx}})_{i,j} = k(\mathbf{x}_i, \mathbf{x}_j), i,j, = 1,2,\cdots,N$, and $\mathbf{I}_N$ is an $N \times N$ identity matrix. With the estimated parameters $\hat{\mu}$, $\hat{\sigma}_Z^2$, $\widehat{\boldsymbol{\theta}}$ and $\hat{\sigma}_n^2$, the predictive mean and variance at a new point $\mathbf{x}^*$ are given by:

$$\hat{\mu}_{\hat{\xi}_j}(\mathbf{x}^*) = \hat{\mu} + \mathbf{k}_{\mathbf{xx}^*}^{\mathrm{T}} (K_{\mathbf{xx}} + \sigma_n^2 \mathbf{I}_N)^{-1} \big(\boldsymbol{\xi}^{(j)} - \mathbf{1}\mu\big), \quad (30)$$

$$\hat{\sigma}_{\hat{\xi}_j}^2(\mathbf{x}^*) = k(\mathbf{x}^*, \mathbf{x}^*) - \mathbf{k}_{\mathbf{xx}^*}^{\mathrm{T}} (K_{\mathbf{xx}} + \sigma_n^2 \mathbf{I}_N)^{-1} \mathbf{k}_{\mathbf{xx}^*}, \quad (31)$$

where $\mathbf{k}_{\mathbf{xx}^*} = [k(\mathbf{x}_1, \mathbf{x}^*), k(\mathbf{x}_2, \mathbf{x}^*), \cdots, k(\mathbf{x}_N, \mathbf{x}^*)]^{\mathrm{T}}$.

The mean and covariance matrix of $\widehat{\boldsymbol{\xi}}(\mathbf{x}^*)$ are $\widehat{\boldsymbol{\mu}}_{\hat{\boldsymbol{\xi}}}(\mathbf{x}^*) = [\hat{\mu}_{\hat{\xi}_1}(\mathbf{x}^*), \cdots, \hat{\mu}_{\hat{\xi}_m}(\mathbf{x}^*)]^{\mathrm{T}}$ and $\widehat{\boldsymbol{\Sigma}}_{\hat{\boldsymbol{\xi}}}(\mathbf{x}^*) = \operatorname{diag}([\hat{\sigma}_{\hat{\xi}_1}^2(\mathbf{x}^*), \cdots, \hat{\sigma}_{\hat{\xi}_m}^2(\mathbf{x}^*)]^{\mathrm{T}})$. And the predicted mean and variance of time-variant response $\hat{y}(t)$ at a specified time node $t^*$ can be obtained as:

$$\hat{\mu}_{\hat{y}}(\mathbf{x}^*, t^*) = \boldsymbol{\eta}(t^*)^{\mathrm{T}} \mathbf{B} \widehat{\boldsymbol{\mu}}_{\hat{\boldsymbol{\xi}}}(\mathbf{x}^*). \quad (32)$$

$$\hat{\sigma}_{\hat{y}}^2(\mathbf{x}^*, t^*) = \boldsymbol{\eta}(t^*)^{\mathrm{T}} \mathbf{B} \widehat{\boldsymbol{\Sigma}}_{\hat{\boldsymbol{\xi}}}(\mathbf{x}^*) \mathbf{B}^{\mathrm{T}} \boldsymbol{\eta}(t^*). \quad (33)$$

Denote $\mathcal{M}_j$ the surrogate model for the $j$-th component of $\boldsymbol{\xi}$, then the surrogate model for the dynamical system can be denoted as $\mathcal{M} = \{\mathcal{M}_1, \mathcal{M}_2, \cdots, \mathcal{M}_m\}$.

### 3.3. Surrogate-based forward and inverse UQ for dynamical systems

Once the surrogate model for the dynamical system is constructed, it can be directly combined with Monte Carlo simulation for forward UQ. Given a set of Monte Carlo samples $\mathbf{x}_1, \mathbf{x}_2, \cdots, \mathbf{x}_{N_{\mathrm{MCS}}}$ according to $f_{\mathbf{X}}(\mathbf{x})$, $\mu_Y(t)$ in Eq. (1) and $\sigma_Y(t)$ in Eq. (2) can be estimated as follows:

$$\hat{\mu}_Y(t) = \frac{1}{N_{\mathrm{MCS}}} \sum_{i=1}^{N_{\mathrm{MCS}}} \boldsymbol{\eta}(t)^{\mathrm{T}} \mathbf{B} \widehat{\boldsymbol{\mu}}_{\hat{\boldsymbol{\xi}}}(\mathbf{x}_i) = \frac{1}{N_{\mathrm{MCS}}} \boldsymbol{\eta}(t)^{\mathrm{T}} \mathbf{B} \sum_{i=1}^{N_{\mathrm{MCS}}} \widehat{\boldsymbol{\mu}}_{\hat{\boldsymbol{\xi}}}(\mathbf{x}_i). \quad (34)$$

$$\hat{\sigma}_Y(t) = \sqrt{\frac{1}{N_{\mathrm{MCS}}} \sum_{i=1}^{N_{\mathrm{MCS}}} \left(\boldsymbol{\eta}(t)^{\mathrm{T}} \mathbf{B} \widehat{\boldsymbol{\mu}}_{\hat{\boldsymbol{\xi}}}(\mathbf{x}_i) - \hat{\mu}_Y(t)\right)^2}, \quad (35)$$

where $\widehat{\boldsymbol{\mu}}_{\hat{\boldsymbol{\xi}}}(\mathbf{x}_i) = [\hat{\mu}_{\hat{\xi}_1}(\mathbf{x}_i), \cdots, \hat{\mu}_{\hat{\xi}_m}(\mathbf{x}_i)]^{\mathrm{T}}$.

For inverse UQ, a Bayesian framework is utilized in this research. The surrogate model $\mathcal{M}$ can be viewed as a function of input parameters **X** mapping to the high-dimensional output **Y**. In the Bayesian framework, a discrepancy term can be added to link predictions $\mathcal{M}(\mathbf{X})$ with observations **Y** as follows:

$$\mathbf{Y} = \mathcal{M}(\mathbf{X}) + \boldsymbol{\epsilon}, \tag{36}$$

where $\boldsymbol{\epsilon}$ is an $N_t$-dimensional vector. For simplicity, we assume $\boldsymbol{\epsilon}$ follows a zero mean multivariate Gaussian distribution with covariance matrix $\sigma^2 \mathbf{I}$, where **I** is an $N_t$-dimensional identity matrix. This assumption makes sense because, for an observed time-variant response, we can consider that zero-mean Gaussian noise with variance $\sigma^2$ is added at each time node due to measurement error. Note that other assumptions about the discrepancy term can also be incorporated into the Bayesian inverse UQ framework.

The posterior distribution of $(\mathbf{X}, \sigma^2)$, under the assumption that **X** and $\sigma^2$ are independent, can be written by Bayes' theorem as:

$$p(\mathbf{X}, \sigma^2 | \mathbf{Y}) \propto p(\mathbf{X}) p(\sigma^2) p(\mathbf{Y} | \mathbf{X}, \sigma^2), \tag{37}$$

where $p(\mathbf{X})$ and $p(\sigma^2)$ are the prior distributions of **X** and $\sigma^2$ respectively, and $p(\mathbf{Y}|\mathbf{X}, \sigma^2)$ is the likelihood function. When $N_{\text{observe}}$ observations are available, $p(\mathbf{Y}|\mathbf{X}, \sigma^2)$ has the following form:

$$p(\mathbf{Y}|\mathbf{X}, \sigma^2) = \prod_{i=1}^{N_{\text{observe}}} \frac{1}{(2\pi\sigma^2)^{N_{\text{observe}}/2}} \exp\left[-\frac{1}{2\sigma^2}(\mathbf{y}_i - \mathcal{M}(\mathbf{X}))^{\mathrm{T}}(\mathbf{y}_i - \mathcal{M}(\mathbf{X}))\right]. \tag{38}$$

In this research, the affine invariant ensemble algorithm [42] is used to calculate the posterior distribution. In the implementation, 100 parallel chains are generated, with initial points randomly sampled from the prior distributions. Each chain is set to run for 300 MCMC iterations. The proportion of samples discarded as burn-in is set to 50%. This MCMC approach is efficiently executed using the UQLab toolbox [43, 44].

## 4. Examples and discussions

In this section, we evaluate the performance of the proposed method alongside several comparative methods on both mathematical and engineering examples. We compare the modeling accuracy of the proposed KFDR with Kriging models that incorporate PCA (KPCA), independent component analysis (KICA), and autoencoders (KAE). Each method follows a similar procedure to KFDR: first performing dimension reduction on the time-variant response, then constructing Kriging models between the inputs and each low-dimensional representation of the response. The primary difference lies in the dimension reduction technique employed by each method. KPCA treats the responses of dynamical systems as vectors and applies standard PCA to reduce response dimensionality. The number of retained principal components is chosen so that they explain more than 99% of the total variance. KICA uses independent component analysis (ICA), a blind source separation technique, for dimension reduction, decomposing a signal into a linear combination of independent component signals. KAE employs an autoencoder for dimension reduction, which is a nonlinear technique with a natural framework for encoding (dimension reduction) and decoding (reconstruction). The transfer functions for both the encoder and decoder are configured as logistic sigmoid functions. The maximum number of training epochs for the autoencoder is set to 1000, and the $L_2$ weight regularization coefficient is set to 0.001. The number of neurons in the hidden layer is set to 20. For the proposed KFDR, we represent the time-variant responses using both Fourier basis functions and B-spline basis functions, referred to as KFDR-F and KFDR-B, respectively. Normalized root mean square error (NRMSE) is used to quantify the modeling error:

$$\text{NRMSE} = \frac{1}{N_{\text{test}}} \sum_{i=1}^{N_{\text{test}}} \left[ \frac{\sqrt{\frac{1}{N_t} \sum_{j=1}^{N_t} \left( y_i(t_j) - \hat{y}_i(t_j) \right)^2}}{\max_j [y_i(t_j)] - \min_j [y_i(t_j)]} \right], \tag{39}$$

where $N_{\text{test}}$ denotes the size of test set, $N_t$ is the number of discretized time nodes, $y_i(t_j)$ and $\hat{y}_i(t_j)$ are the true and predicted values of the $i$-th time-variant response at $t_j$, respectively. Besides, the performances of different methods on forward and inverse uncertainty quantification tasks are investigated.

### 4.1. Example 1: The Duffing oscillator

The Duffing oscillator adopted from [17] is used as the first example. The governing ordinary differential equation for the Duffing oscillator is as follows:

$$m\ddot{y}(t) + c\dot{y}(t) + ky(t) + k_2 y^2(t) + k_3 y^3(t) = f(t), \tag{40}$$

where $m = 1$, $k = 1 \times 10^4$, $k_2 = 1 \times 10^7$, $k_3 = 5 \times 10^9$, $y(t)$ is the displacement of oscillator with initial conditions $\dot{y}(0) = 0$ and $y(0) = y_0$, and $f(t)$ is the excitation force given by:

$$f(t) = \alpha \cos(\beta t) + \sin\bigl((\beta + 3)t\bigr) + \sin(2\beta t). \tag{41}$$

Physical units have been dropped intentionally for simplicity. The quantity of interest is the oscillator displacement $y(t)$ over the time interval [0, 2]. Runge-Kutta method is used to solve Eq. (40) to obtain $y(t)$, and the time interval is uniformly discretized into 401 time nodes. The parameters $\alpha$, $\beta$, $c$, and $y_0$ are set as input variables, with their lower and upper bounds listed in Table 4. Fig. 4 shows 100 different realizations for this problem.

**Table 4**
Lower and upper bounds of inputs of the Duffing oscillator.

| Variables | Lower bounds | Upper bounds |
|---|---|---|
| $\alpha$ | 0.6 | 1.4 |
| $\beta$ | 1.5 | 2.5 |
| $c$ | 0.6 | 1.4 |
| $y_0$ | $-1 \times 10^{-4}$ | 0 |

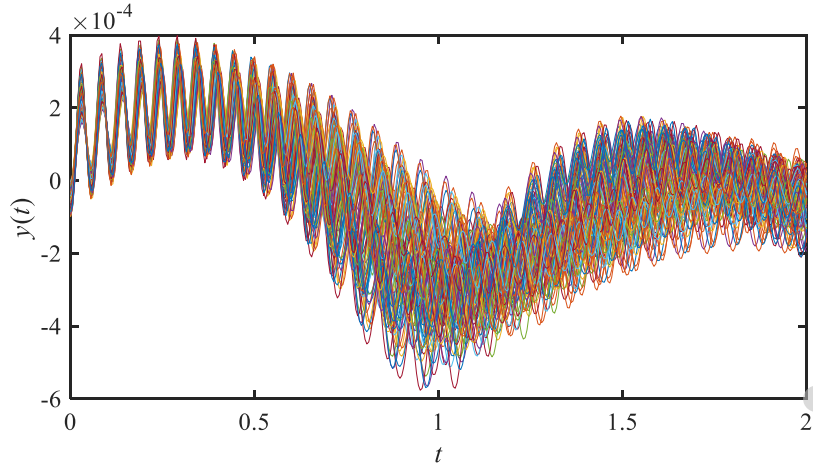

**Fig. 4**. 100 realizations of the responses for the Duffing oscillator problem.

The number of basis functions ($N_b$) for KFDR-F and KFDR-B is 201 and 405, respectively, while the number of retained latent functions ($m$) for both KFDR-F and KFDR-B is 14. Fig. 5 illustrates the modeling error of different methods across different training sample sizes. In the left panel of Fig. 5, each box displays the median as the central mark, with the bottom and top edges representing the 25th and 75th percentiles, respectively. The solid lines extend to the most extreme data points that are not considered outliers, while outliers are indicated separately using diamond markers. The error is evaluated using an additional test set of 1000 samples, with both training and test samples generated through Latin hypercube sampling. To mitigate randomness effects, each experiment is repeated ten times. In each trial, identical training samples are used to construct surrogate models for all methods.

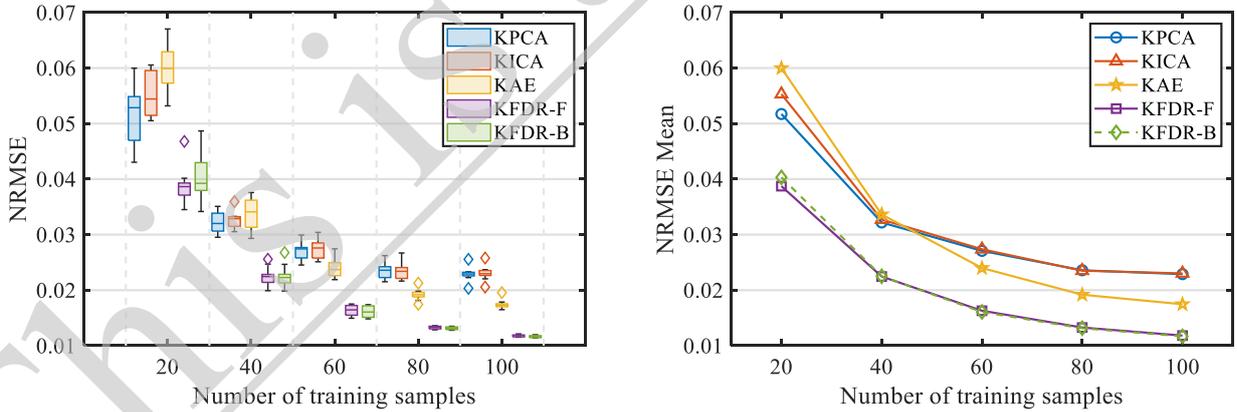

**Fig. 5**. Boxplots (left) and means (right) of the normalized root mean square errors of different methods across different training sample sizes for the Duffing oscillator problem.

Fig. 5 demonstrates a clear downward trend in NRMSE across all methods as the number of training samples increases. However, the proposed KFDR-F and KFDR-B yield a smaller NRMSE compared to KPCA, KICA, and KAE across all training sample sizes, indicating higher modeling accuracy. The proposed approach outperforms KPCA and KICA because KPCA and KICA are based on linear dimensionality reduction techniques, while KFDR-F and KFDR-B can capture nonlinear features in the response through basis expansion in the functional space. As a result, the proposed method offers a more flexible representation than the linear methods. Although the autoencoder is a powerful nonlinear dimensionality reduction method, it may lose effectiveness with a small sample size. Consequently, the modeling accuracy of KAE is not as high as that of KFDR-F and KFDR-B. An interesting phenomenon is that as the number of training samples increases, the modeling error of KAE becomes smaller than that of KPCA and KICA. This is because the autoencoder can extract nonlinear features more effectively with a large sample size, highlighting the potential of neural network-based approaches when handling large datasets.

In addition, we investigate the influence of noise level $\sigma$ and training sample size on the modeling accuracy of the proposed KFDR method. Zero-mean Gaussian noise with varying standard deviations ($\sigma$) is added to the training output data. The results based on KFDR-B are presented in Fig. 6, which depicts the NRMSE as a function of the training set size for $\sigma = 1\times10^{-5}$, $5\times10^{-5}$, and $1\times10^{-4}$. Additionally, we compare the proposed method to the approach that does not include the roughness regularization term in Eq. (13). It is observed that NRMSE decreases as the number of training samples increases across all noise levels and methods. For all training sample size, larger $\sigma$ values result in higher NRMSE, indicating the increased challenge of accurate modeling under noisy conditions. Methods with regularization (solid lines) exhibit consistently lower NRMSE compared to those without regularization (dashed lines), demonstrating the effectiveness of the roughness regularization term in enhancing model robustness, particularly in noisy scenarios. Furthermore, for smaller $\sigma$ values ($1\times10^{-5}$), the performance gap between methods with and without regularization is less significant. However, at higher noise levels ($1\times10^{-4}$), the benefit of regularization becomes more evident, highlighting its importance in handling noisy data effectively.

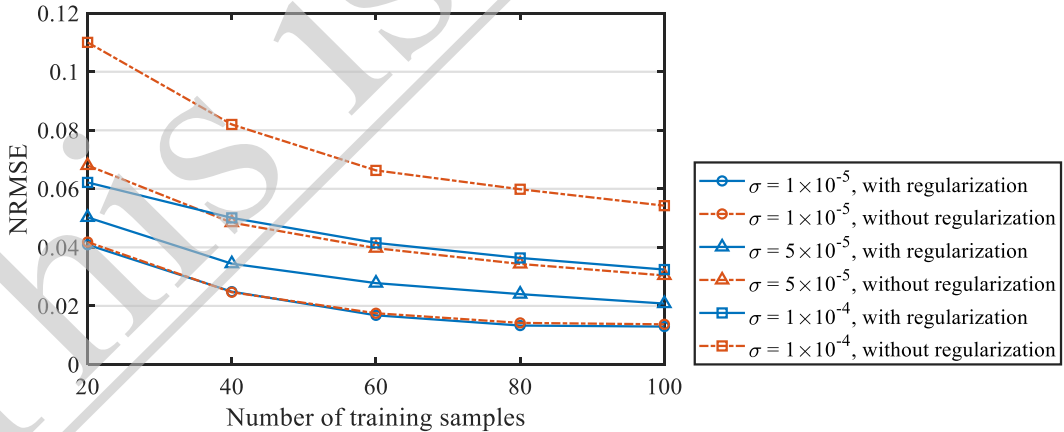

**Fig. 6**. Normalized root mean square errors for different noise levels with and without regularization as a function of the number of training samples for the Duffing oscillator problem.

For the forward uncertainty quantification, the uncertainty information of the input parameters is provided in Table 5. Forward UQ is conducted using the real model and surrogate models trained on 100 samples with different methods. The number of Monte Carlo simulation samples for forward UQ is $1\times10^5$.

Since the modeling accuracy of KPCA and KICA are close, only KPCA is used for forward UQ. Likewise, only KFDR-B is selected for forward UQ between KFDR-F and KFDR-B. Fig. 7 illustrates the forward UQ results. From the upper left panel, we can see that all methods provide accurate predictions of the mean function of the dynamical system's response. While KFDR-B can obtain a more accurate estimation of the standard deviation function than other methods. The lower two panels of Fig. 6 show the probability density functions of the maximum and minimum time-variant responses, fitted using the kernel density estimation method. The probability density function obtained by KFDR-B is closer to the true probability density function than those obtained by other methods, indicating that the proposed approach can achieve higher accuracy in the forward UQ task.

**Table 5**
Uncertainty information of the parameters of the Duffing oscillator.

| Variables | Distribution | Mean | Standard deviation |
|---|---|---|---|
| $\alpha$ | Normal | 1.0 | 0.05 |
| $\beta$ | Normal | 2.0 | 0.1 |
| $c$ | Normal | 1.0 | 0.05 |
| $y_0$ | Normal | $-5\times10^{-5}$ | $5\times10^{-6}$ |

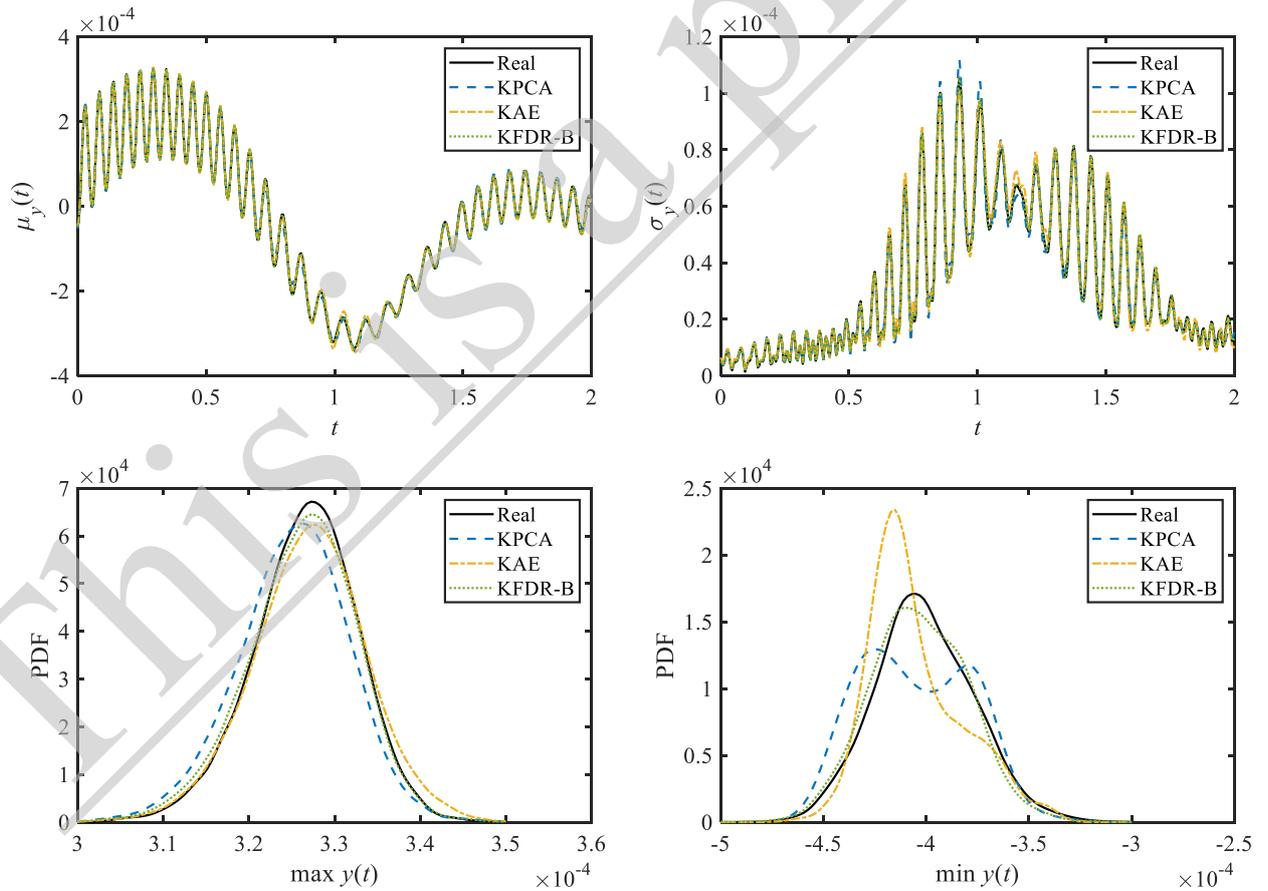

**Fig. 7**. Mean functions over time (upper left), standard deviation functions over time (upper right), maximum value distributions (lower left), and minimum value distributions (lower right) of real and predicted time-variant responses for the Duffing oscillator problem.

For inverse uncertainty quantification, the four parameters $\alpha$, $\beta$, $c$, and $y_0$ are assumed to follow uniform prior distributions, with their lower and upper bounds provided in Table 4. The data for inverse UQ consists of three observations at $[\alpha, \beta, c, y_0] = [1.19, 1.82, 0.94, -3.3 \times 10^{-5}]$, with zero-mean Gaussian noise having a standard deviation of $1 \times 10^{-5}$ added at each time node. Table 6 presents the inverse UQ results, showing the mean values and 95% credible intervals of the calibration parameters. Fig. 8 shows the posterior distributions of the calibration parameters. The results indicate that the posterior distributions of $c$ and $y_0$ obtained using the KPCA method exhibit a significant deviation from those of the real model. Similarly, the posterior distribution of $c$ obtained using the KAE method shows a notable deviation. In contrast, the posterior distributions obtained using the KFDR-B method are closer to those of the real model than those from KPCA and KAE, effectively inferring the correct distributions of the calibration parameters. This demonstrates that the proposed approach can achieve higher accuracy in the inverse UQ task.

**Table 6**
Inverse uncertainty quantification results of the Duffing oscillator.

| Variables | Methods | Mean values | 95% credible intervals |
|---|---|---|---|
| $\alpha$ | Real | 1.1927 | [1.1894, 1.1958] |
| | KPCA | 1.1850 | [1.1782, 1.1916] |
| | KAE | 1.1964 | [1.1917, 1.2011] |
| | KFDR-B | 1.1878 | [1.1847, 1.1910] |
| $\beta$ | Real | 1.8200 | [1.8190, 1.8210] |
| | KPCA | 1.8172 | [1.8152, 1.8192] |
| | KAE | 1.8140 | [1.8119, 1.8163] |
| | KFDR-B | 1.8191 | [1.8180, 1.8202] |
| $c$ | Real | $9.4627 \times 10^{-1}$ | $[9.1448, 9.8181] \times 10^{-1}$ |
| | KPCA | $7.5348 \times 10^{-1}$ | $[6.8398, 8.1570] \times 10^{-1}$ |
| | KAE | $7.7070 \times 10^{-1}$ | $[7.0421, 8.5689] \times 10^{-1}$ |
| | KFDR-B | $9.8347 \times 10^{-1}$ | $[9.4824, 10.179] \times 10^{-1}$ |
| $y_0$ | Real | $-3.2712 \times 10^{-5}$ | $[-3.4977, -3.0616] \times 10^{-5}$ |
| | KPCA | $-2.2209 \times 10^{-5}$ | $[-2.5622, -1.8702] \times 10^{-5}$ |
| | KAE | $-3.0564 \times 10^{-5}$ | $[-3.3123, -2.8259] \times 10^{-5}$ |
| | KFDR-B | $-3.3765 \times 10^{-5}$ | $[-3.5846, -3.1583] \times 10^{-5}$ |

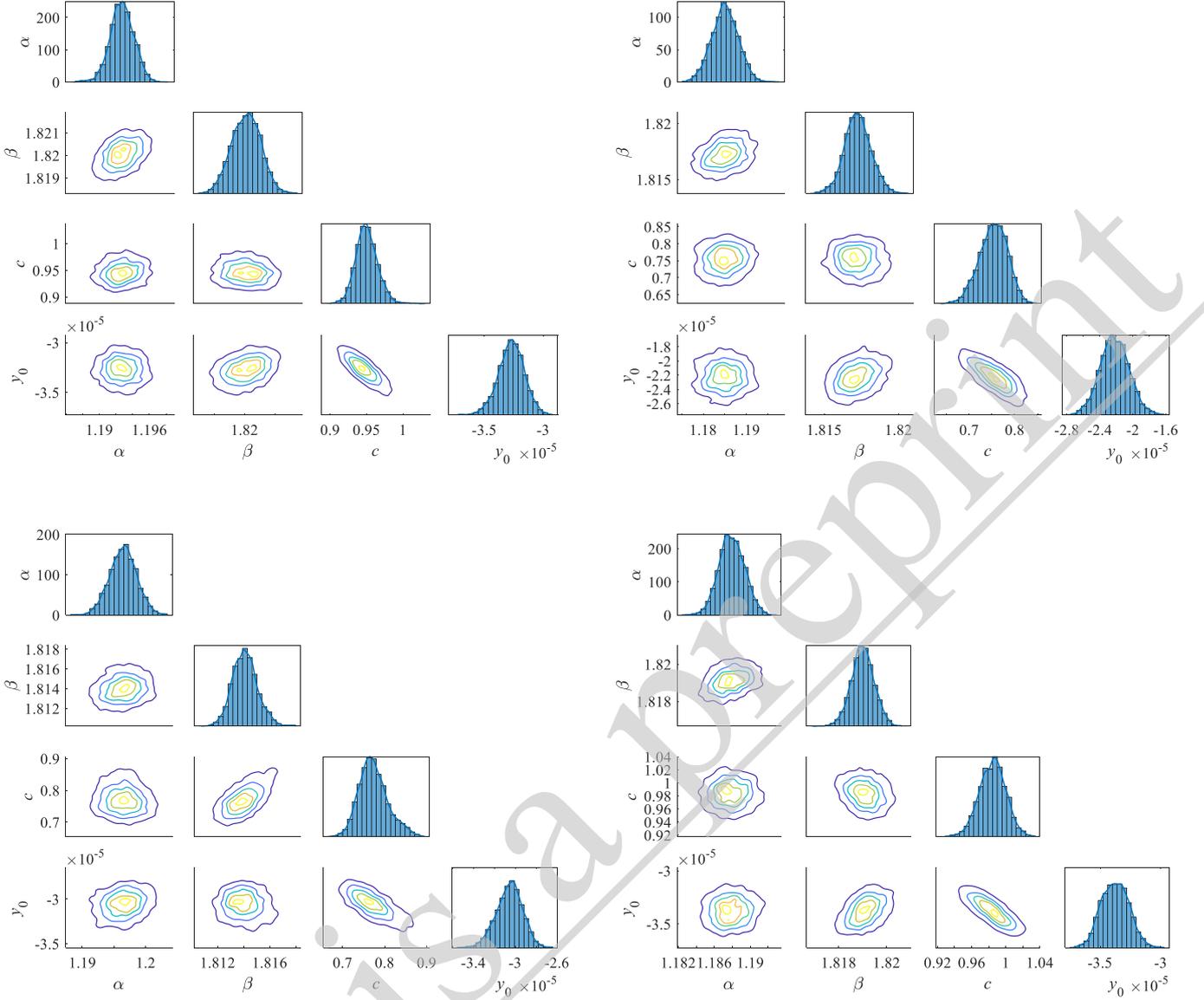

**Fig. 8**. Posterior distributions of the four calibration parameters for the Duffing oscillator problem: real model (upper left), KPCA model (upper right), KAE model (lower left), and proposed KFDR-B model (lower right).

### 4.2. Example 2: The Bouc-Wen hysteretic oscillator

In this example, the forward and inverse UQ of a nonlinear Bouc-Wen oscillator [15] are investigated. The Bouc-Wen model is described by the following differential equation:

$$\begin{cases} m\ddot{y}(t) + c\dot{y}(t) + k[\alpha y(t) + (1-\alpha)z(t)] = f(t), \\ \dot{z}(t) = A\dot{y}(t) - \beta|\dot{y}(t)||z(t)|^{n-1}z(t) - \gamma\dot{y}(t)|z(t)|^n, \end{cases} \quad (42)$$

where $m$ is the mass, $y(t)$ is the displacement of oscillator with initial conditions $\dot{y}(0) = 0$ and $y(0) = y_0$, $c$ the viscous damping coefficient, $k$ the stiffness, $\alpha$ the degree of hysteresis, $z(t)$ the hysteretic displacement with zero initial condition, $f(t)$ the excitation force, and $A, \beta, \gamma, n$ are parameters controlling the behavior of hysteresis and are set $A = 1$, $\beta = \gamma = 7.8 \times 10^3$, $n = 3$. In this example, the excitation force is fixed in the following form:

$$f(t) = -\sqrt{0.006\pi}m \sum_{k=1}^{150}[\vartheta_k \cos(0.1\pi kt) + \vartheta_{150+k} \sin(0.1\pi kt)], \tag{43}$$

where $\vartheta_k$ is a realization of the standard normal distribution. The quantity of interest is the oscillator displacement $y(t)$ over the time interval [0, 16]. Runge-Kutta method is used to solve Eq. (42) to obtain $y(t)$, and the time interval is uniformly discretized into 401 time nodes. The parameters $m$, $c$, $k$, $\alpha$ and $y_0$ are set as input variables, with their lower and upper bounds listed in Table 7. Fig. 9 shows 100 realizations of the responses for this problem.

**Table 7**
Lower and upper bounds of inputs of the Bouc-Wen oscillator.

| Variables | Lower bounds | Upper bounds |
|---|---|---|
| $m$ (kg) | $4 \times 10^4$ | $8 \times 10^4$ |
| $c$ (kg/s) | $8 \times 10^4$ | $1.2 \times 10^5$ |
| $k$ (N/m) | $4 \times 10^6$ | $6 \times 10^6$ |
| $\alpha$ | 0.1 | 0.3 |
| $y_0$ (m) | -0.02 | 0.02 |

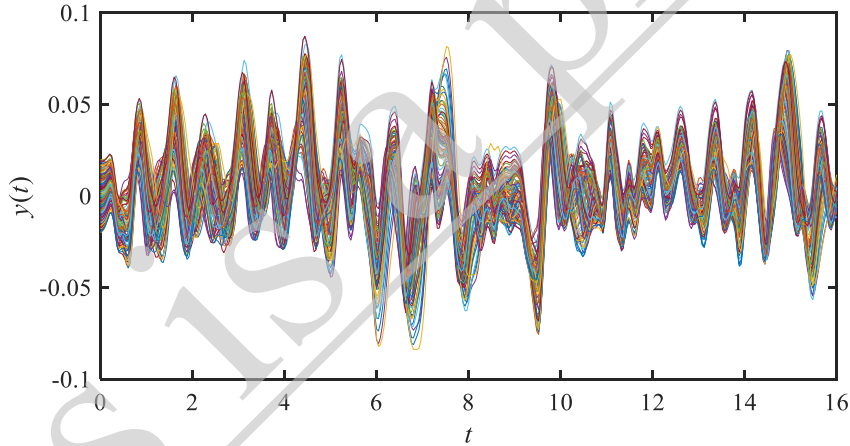

**Fig. 9**. 100 realizations of the responses for the Bouc-Wen oscillator problem.

The number of basis functions ($N_b$) for KFDR-F and KFDR-B is 201 and 405, respectively, while the number of retained latent functions ($m$) for both KFDR-F and KFDR-B is 7. Fig. 10 shows the modeling error of various methods across different training sample sizes, evaluated on a test set of 1000 samples generated with Latin hypercube sampling. Each experiment is repeated ten times to reduce randomness effects. Fig. 10 demonstrates the proposed KFDR-B yield a smaller NRMSE compared to other methods across all training sample sizes, indicating higher modeling accuracy. However, KFDR-F performs poorly in this example, likely because Fourier basis systems are not well-suited for capturing the motion of the Bouc-Wen oscillator. This indicates that the B-spline basis system is more flexible than Fourier basis systems and can represent a broader range of functions. Again, as the number of training samples increases, the modeling error of KAE becomes smaller than that of KPCA and KICA but remains larger than that of KFDR-B, demonstrating the advantage of the proposed approach when faced with a small dataset.

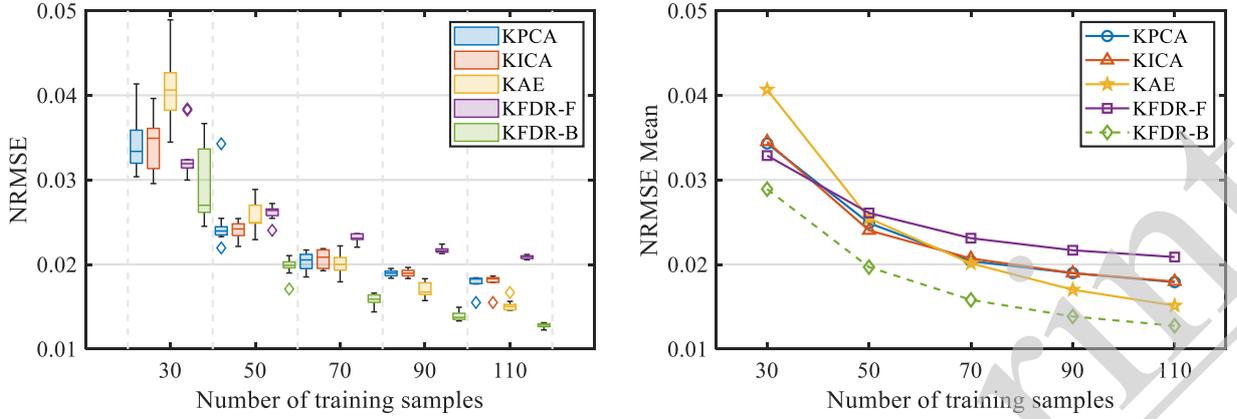

**Fig. 10**. Boxplots (left) and means (right) of the normalized root mean square errors of different methods across different training sample sizes for the Bouc-Wen oscillator problem.

In addition, we investigate the influence of noise level $\sigma$ and training sample size on the modeling accuracy of the proposed KFDR method. Zero-mean Gaussian noise with varying standard deviations ($\sigma = 1\times10^{-3}$, $5\times10^{-3}$, and $1\times10^{-2}$) is added to the training output data. The results based on KFDR-B are presented in Fig. 11, which depicts the NRMSE as a function of the training sample size for different $\sigma$. Additionally, we compare the proposed method to the approach that does not include the roughness regularization. It is observed that NRMSE decreases as the number of training samples increases across all noise levels and methods. For all training sample size, larger $\sigma$ values result in higher NRMSE, indicating the increased challenge of accurate modeling under noisy conditions. For smaller $\sigma$ values ($1\times10^{-3}$), the performance gap between methods with and without regularization is less significant. However, at higher noise levels ($1\times10^{-2}$), the benefit of regularization becomes more evident, demonstrating the effectiveness of the roughness regularization term in enhancing model robustness, particularly in large noisy scenarios.

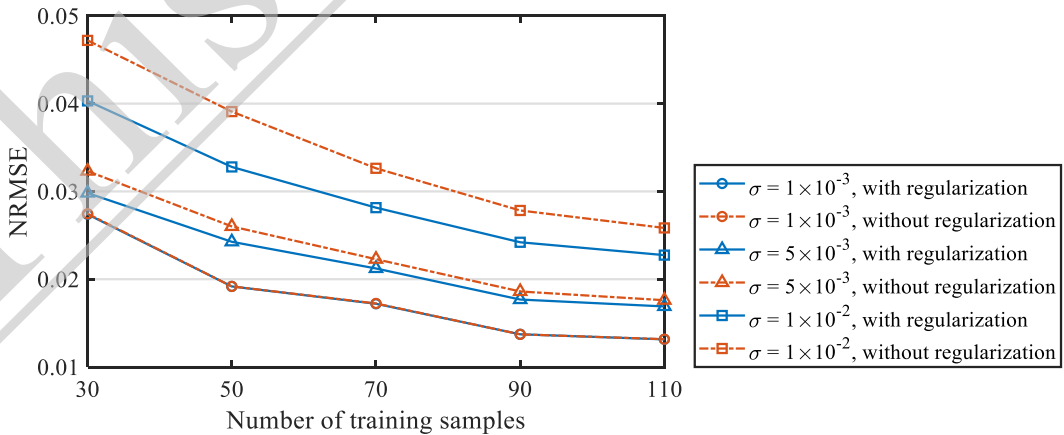

**Fig. 11**. Normalized root mean square errors for different noise levels with and without regularization as a function of the number of training samples for the Bouc-Wen oscillator problem.

For the forward UQ, the uncertainty information of the input parameters is provided in Table 8. Forward UQ is conducted using the real model and surrogate models trained on 110 samples with KPCA, KAE, and KFDR-B. The number of Monte Carlo simulation samples for forward UQ is $1\times10^5$. Fig. 12 shows the forward UQ results. Again, all methods provide accurate predictions of the mean function of the dynamical system's response. While KFDR-B obtains a more accurate estimation of the standard deviation function than other methods. The lower two panels of Fig. 12 show that the extreme value distributions obtained by KFDR-B are closer to the true probability density function than those from other methods, indicating that the proposed approach achieves higher accuracy in the forward UQ task.

**Table 8**
Uncertainty information of the parameters of the Bouc-Wen oscillator.

| Variables | Distribution | Mean | Standard deviation |
|---|---|---|---|
| $m$ (kg) | Lognormal | $6\times10^4$ | $3\times10^3$ |
| $c$ (kg/s) | Lognormal | $1\times10^5$ | $3\times10^3$ |
| $k$ (N/m) | Lognormal | $5\times10^6$ | $1\times10^5$ |
| $\alpha$ | Normal | 0.2 | 0.01 |
| $y_0$ (m) | Normal | 0 | 0.002 |

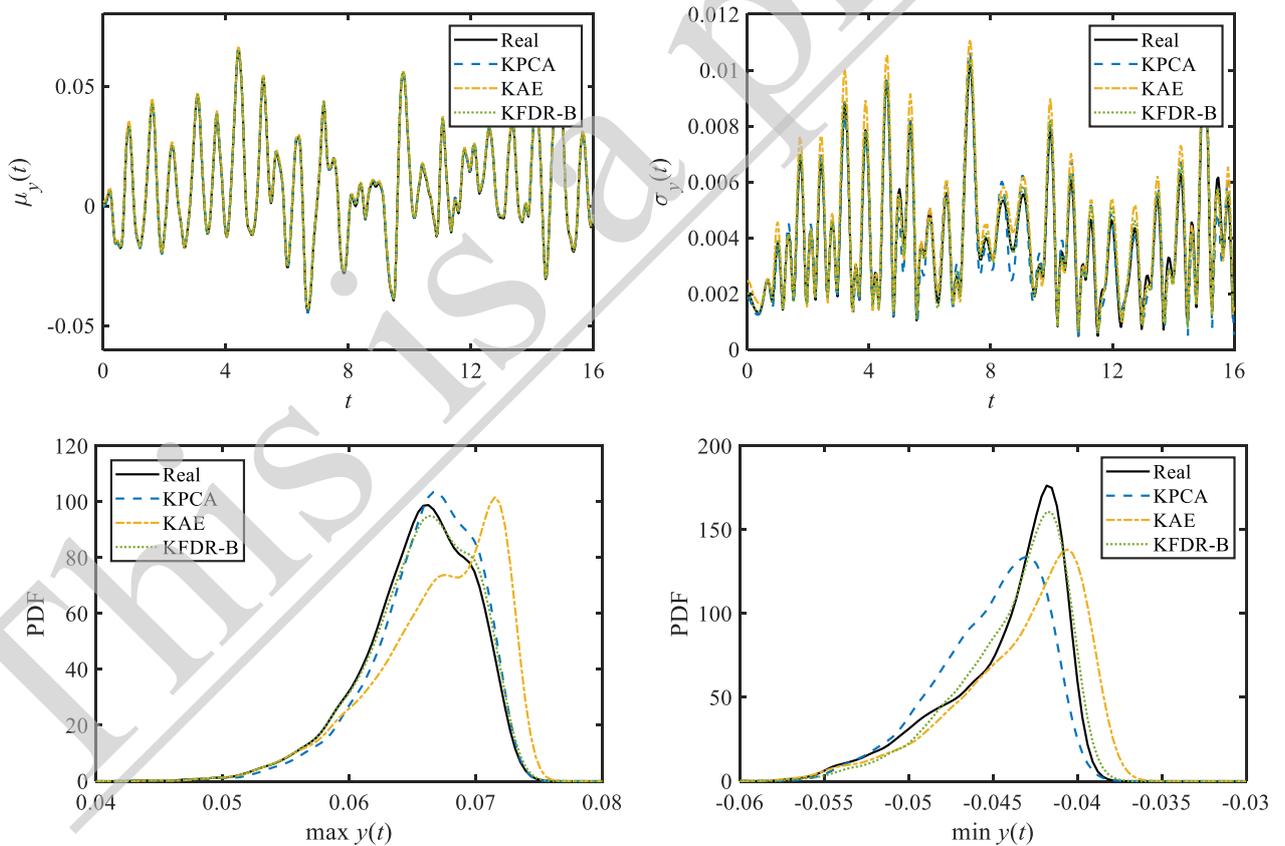

**Fig. 12**. Mean functions over time (upper left), standard deviation functions over time (upper right), maximum value distributions (lower left), and minimum value distributions (lower right) of real and predicted time-variant responses for the Bouc-Wen oscillator problem.

For inverse uncertainty quantification, the mass of oscillator is fixed at $7\times10^4$, the other four parameters $c$, $k$, $\alpha$ and $y_0$ are assumed to follow uniform prior distributions, with their lower and upper bounds provided in Table 7. The data for inverse UQ consists of three observations at $[c, k, \alpha, y_0] = [1.05 \times 10^5, 4.77 \times 10^6, 0.21, 0.01]$, with zero-mean Gaussian noise having a standard deviation of $5\times10^{-3}$ added at each time node. Table 9 presents the inverse UQ results, showing the mean values and 95% credible intervals of the calibration parameters. Fig. 13 shows the posterior distributions of the calibration parameters. The results indicate that all methods provide relatively accurate posterior distributions for $c$ and $k$. However, the posterior distributions of $\alpha$ obtained using the KPCA and KAE methods show a significant deviation from those of the real model. Additionally, KPCA and KAE produce wider 95% credible intervals for $y_0$ compared to the real model and the KFDR-B method. Moreover, methods KPCA and KAE erroneously infer a strong positive correlation between $\alpha$ and $y_0$. In contrast, the KFDR-B method yields posterior distributions for $\alpha$ and $y_0$ that are very close to those of the real model, once again demonstrating the high accuracy of the proposed method in inverse UQ.

**Table 9**
Inverse uncertainty quantification results of the Bouc-Wen oscillator.

| Variables | Methods | Mean values | 95% credible intervals |
|---|---|---|---|
| $c$ | Real | $1.0679 \times 10^5$ | $[1.0415, 1.0922] \times 10^5$ |
| | KPCA | $1.0999 \times 10^5$ | $[1.0661, 1.1325] \times 10^5$ |
| | KAE | $1.0968 \times 10^5$ | $[1.0550, 1.1349] \times 10^5$ |
| | KFDR-B | $1.0722 \times 10^5$ | $[1.0468, 1.0971] \times 10^5$ |
| $k$ | Real | $4.7655 \times 10^6$ | $[4.7510, 4.7812] \times 10^6$ |
| | KPCA | $4.7438 \times 10^6$ | $[4.7254, 4.7614] \times 10^6$ |
| | KAE | $4.7529 \times 10^6$ | $[4.7324, 4.7741] \times 10^6$ |
| | KFDR-B | $4.7638 \times 10^6$ | $[4.7487, 4.7788] \times 10^6$ |
| $\alpha$ | Real | $2.1378 \times 10^{-1}$ | $[2.0790, 2.1993] \times 10^{-1}$ |
| | KPCA | $2.2865 \times 10^{-1}$ | $[1.9450, 2.6527] \times 10^{-1}$ |
| | KAE | $2.2429 \times 10^{-1}$ | $[2.1517, 2.3361] \times 10^{-1}$ |
| | KFDR-B | $2.1250 \times 10^{-1}$ | $[2.0689, 2.1788] \times 10^{-1}$ |
| $y_0$ | Real | $0.9854 \times 10^{-2}$ | $[0.9001, 1.0740] \times 10^{-2}$ |
| | KPCA | $1.1253 \times 10^{-2}$ | $[0.7830, 1.5239] \times 10^{-2}$ |
| | KAE | $1.0353 \times 10^{-2}$ | $[0.8838, 1.2022] \times 10^{-2}$ |
| | KFDR-B | $0.9757 \times 10^{-2}$ | $[0.8903, 1.0618] \times 10^{-2}$ |

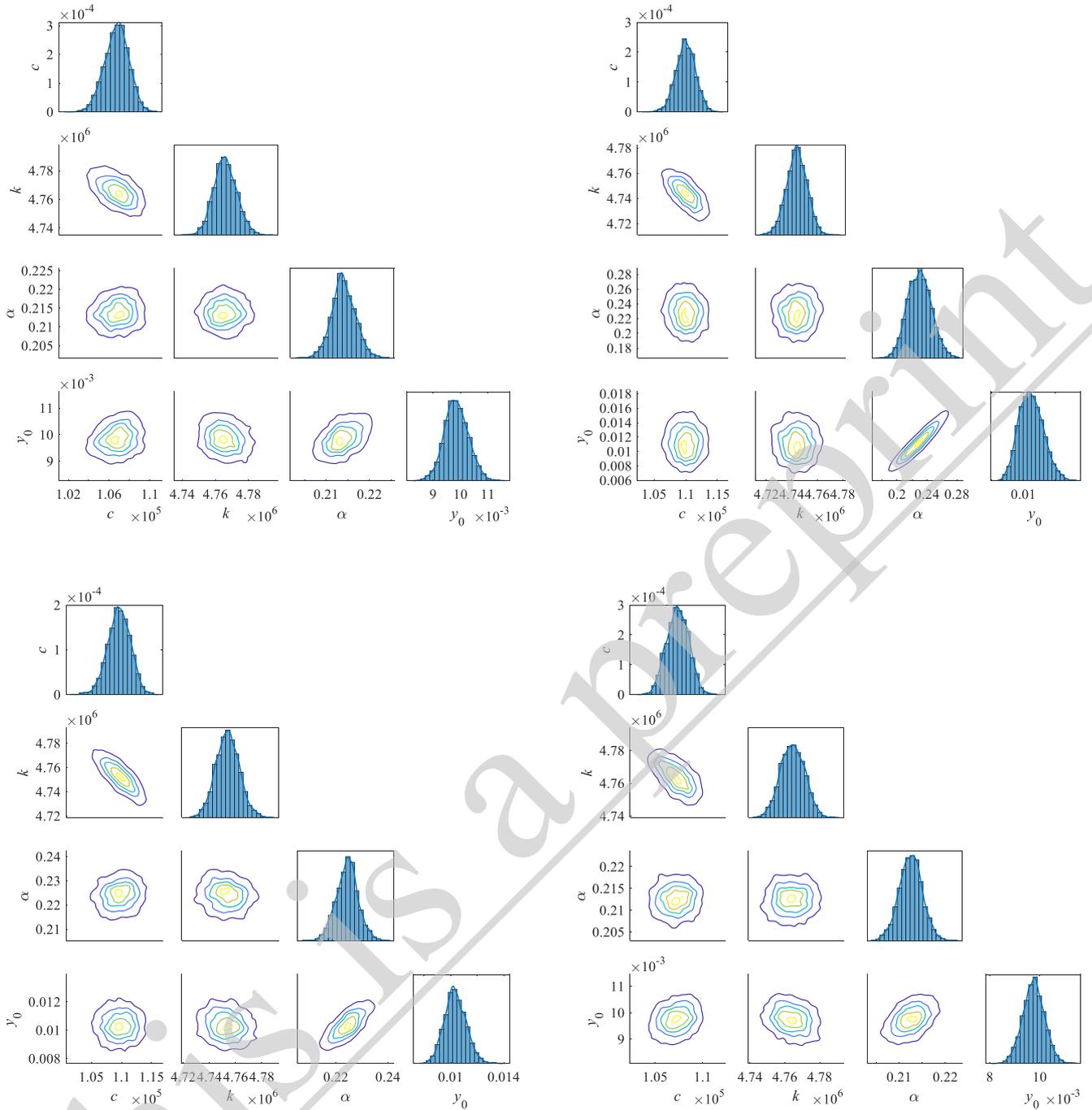

**Fig. 13**. Posterior distributions of the four calibration parameters for the Bouc-Wen oscillator problem: real model (upper left), KPCA model (upper right), KAE model (lower left), and proposed KFDR-B model (lower right).

## 4.3. Example 3: A crane structure

This example considers the transient analysis of a crane structure under a sudden load. Fig. 14 shows a schematic and Fig 15 shows the dimensions of the crane. The crane is composed of steel box beams with two different cross-sections for the main beams and bracing beams. One end of the main beams is fixed at points A, B, C, and D, while the other end (point E) is subjected to an instantaneous impact force with a magnitude of $F$ and a duration of $T_F$. All beams are made of the same steel material, with density $\rho$, Young's modulus $E$, and shear modulus $G$ treated as varying parameters. In addition, $F$ and $T_F$ are also treated as varying parameters. Table 10 presents the lower and upper bounds for these five input parameters. The quantity of interest is the force in the Y-direction at point A over the specified time interval [0, 0.5s], which is obtained through finite element analysis (FEA). The time interval is uniformly discretized into 201 time nodes.

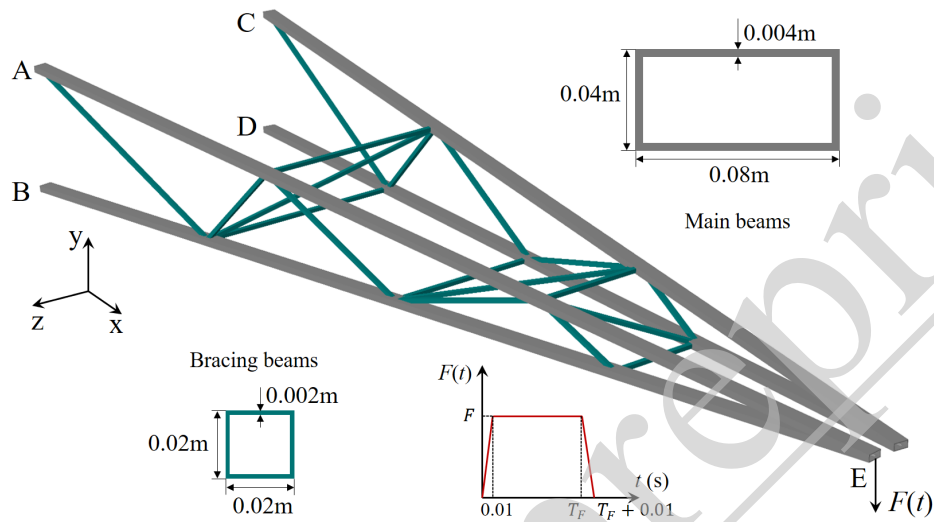

**Fig. 14**. The crane structure subjected to an instantaneous impact force.

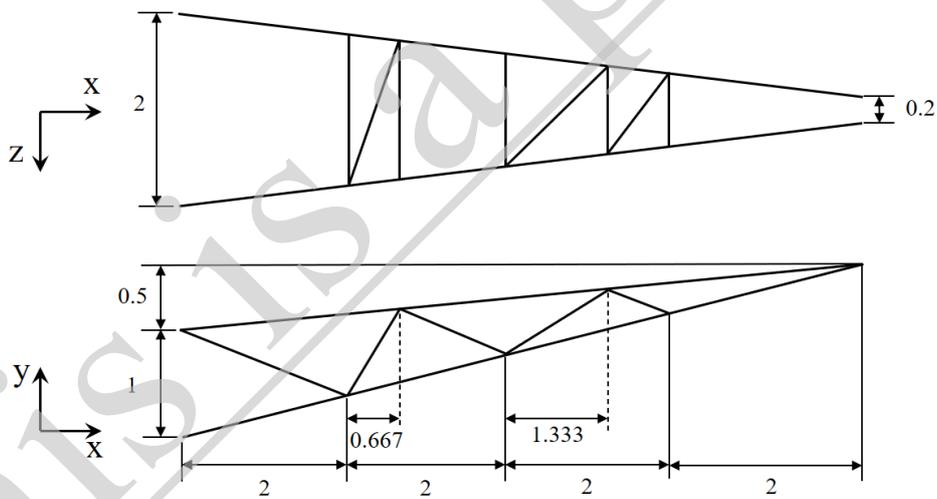

**Fig. 15**. Dimensions of the crane structure (unit: meters).

**Table 10**
Lower and upper bounds of inputs of the crane structure.

| Variables | Lower bounds | Upper bounds |
|---|---|---|
| $\rho$ (kg/m³) | 7600 | 8000 |
| $E$ (Pa) | $1.8\times10^{11}$ | $2.2\times10^{11}$ |
| $G$ (Pa) | $7.8\times10^{10}$ | $8.2\times10^{10}$ |
| $F$ (N) | $-12\times10^3$ | $-6\times10^3$ |
| $T_F$ (s) | 0.17 | 0.23 |

We collected 100 samples using FEA, with input samples generated through Latin hypercube sampling. Fig. 16 illustrates the 100 realizations of responses. Ten-fold cross-validation was employed to evaluate the modeling accuracy of the various methods in this example, and the process was repeated ten times to mitigate the impact of randomness. The number of basis functions ($N_b$) for KFDR-F and KFDR-B is 151 and 205, respectively, while the number of retained latent functions ($m$) for both KFDR-F and KFDR-B is 11. Fig. 17 illustrates the modeling errors of the different methods. The proposed methods, KFDR-F and KFDR-B, show lower NRMSE values compared to the comparative methods (KPCA, KICA, and KAE), indicating better modeling accuracy. Additionally, their narrower boxplots suggest more consistent performance across different trials. Among the comparative methods, KPCA and KICA have higher NRMSE values with larger variability, reflecting lower accuracy and stability. While KAE achieves lower median NRMSE than KPCA and KICA.

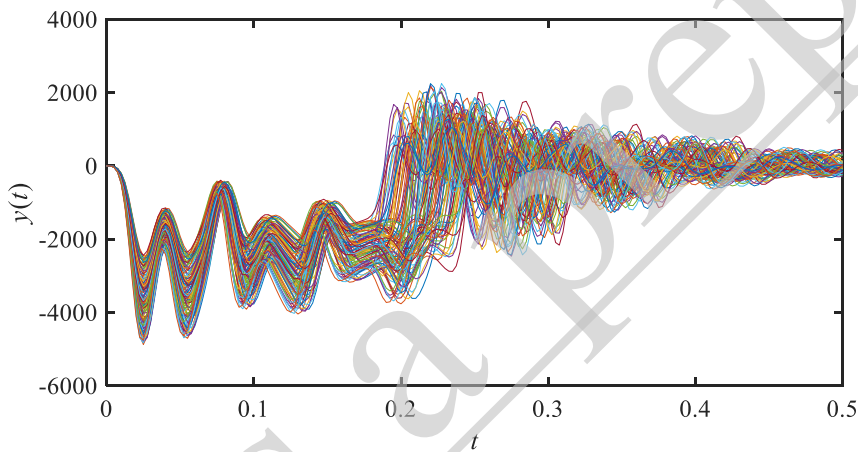

**Fig. 16**. 100 realizations of the responses for the crane structure problem.

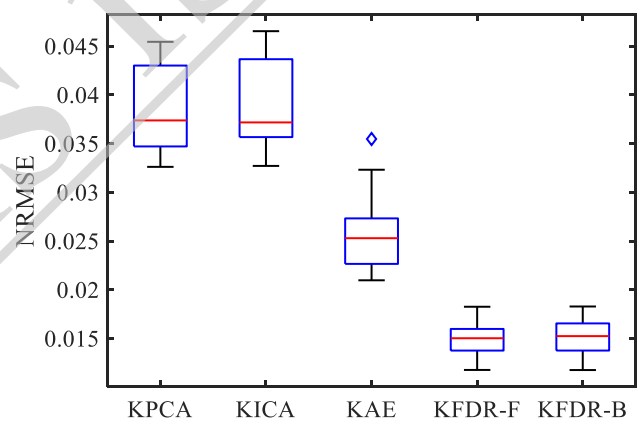

**Fig. 17**. Boxplots of the normalized root mean square errors of different methods for the crane structure problem.

For the forward UQ, the uncertainty information of the input parameters is provided in Table 11. Forward UQ is conducted using the surrogate models trained on all 100 samples with KPCA, KAE, and KFDR-B.

The number of Monte Carlo simulation samples for forward UQ is $1\times10^5$. Fig. 18 presents the forward UQ results, showing that the mean functions predicted by the three methods are consistent, while the standard deviation functions exhibit differences among the methods. All standard deviation functions exhibit higher values around 0.2s, as the external force is removed at this point, causing the crane to transition from forced vibration to free vibration. For the extreme value distributions, the three methods predict different modes for the maximum value distribution, with KAE even producing a multimodal PDF. In contrast, all three methods predict the same mode for the minimum value distribution, although the PDF obtained by KAE differs from those of KPCA and KFDR-B.

**Table 11**
Uncertainty information of the parameters of the crane structure.

| Variables | Distribution | Mean | Standard deviation |
| --- | --- | --- | --- |
| $\rho$ (kg/m$^3$) | Lognormal | 7800 | 20 |
| $E$ (Pa) | Lognormal | $2\times10^{11}$ | $2.5\times10^9$ |
| $G$ (Pa) | Lognormal | $8\times10^{10}$ | $2\times10^8$ |
| $F$ (N) | Normal | $-9\times10^3$ | 500 |
| $T_F$ (s) | Lognormal | 0.2 | 0.005 |

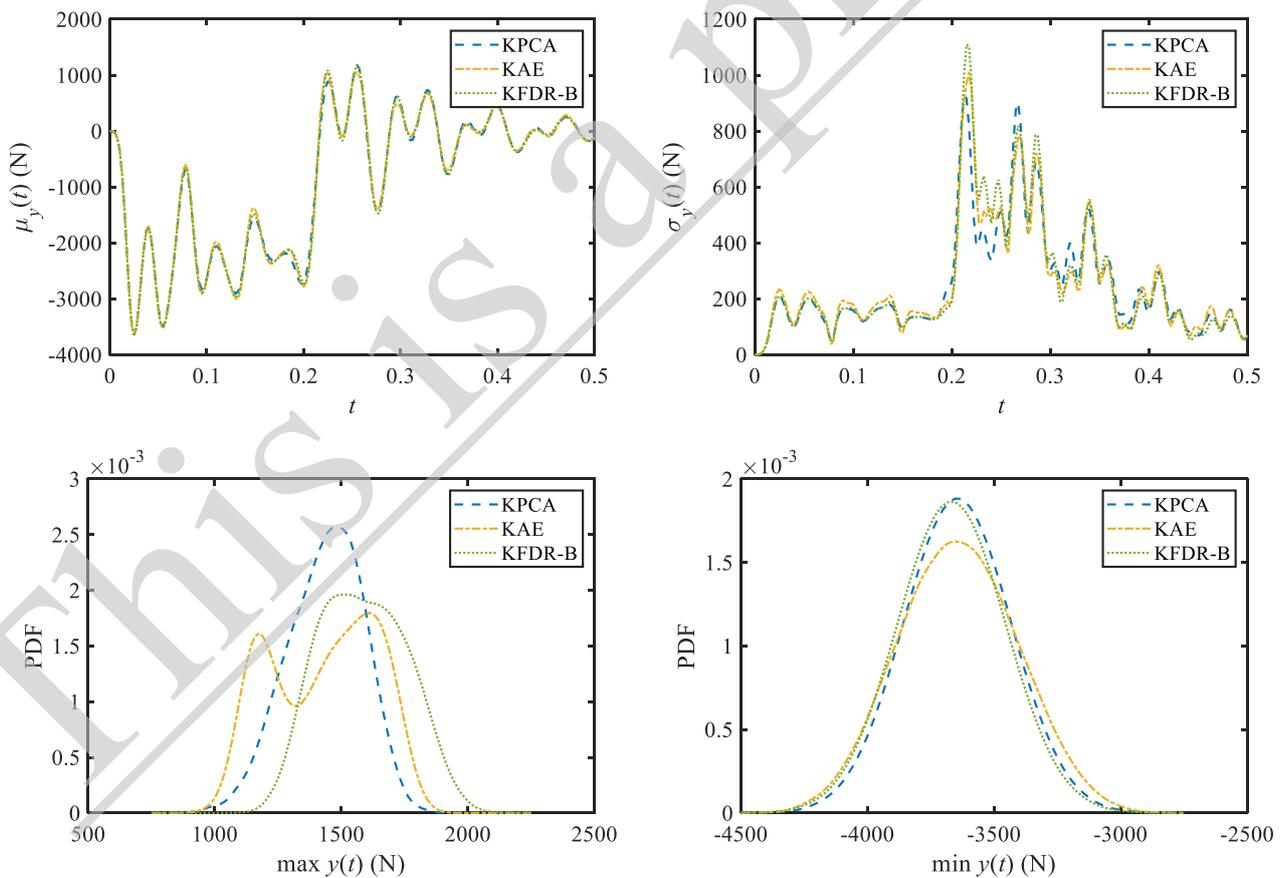

**Fig. 18**. Mean functions over time (upper left), standard deviation functions over time (upper right), maximum value distributions (lower left), and minimum value distributions (lower right) of predicted time-variant responses for the crane structure problem.

For inverse uncertainty quantification, $\rho$, $E$, and $G$ are fixed at 7800, $2\times10^{11}$, and $8\times10^{10}$, respectively. $F$ and $T_F$ are assumed to follow uniform prior distributions in $[-11\times10^3, -9\times10^3]$ and $[0.19, 0.23]$, respectively. The data for inverse UQ consists of three observations at $[F, T_F] = [-9.8\times10^3, 0.21]$, with zero-mean Gaussian noise having a standard deviation of 100 added at each time node. Table 12 presents the inverse UQ results, showing the mean values and 95% credible intervals of the calibration parameters. Fig. 19 shows the posterior distributions of the calibration parameters. The results show that the proposed KFDR-B method generates posterior distributions that are very close to the true values, whereas the KPCA and KAE methods exhibit slight deviations from the true values.

**Table 12**
Inverse uncertainty quantification results of the Bouc-Wen oscillator.

| Variables | Methods | Mean values | 95% credible intervals |
|---|---|---|---|
| $F$ | KPCA | $-9.7510\times10^3$ | $[-9.7997, -9.7036]\times10^3$ |
|  | KAE | $-9.8694\times10^3$ | $[-9.9484, -9.7885]\times10^3$ |
|  | KFDR-B | $-9.8083\times10^3$ | $[-9.8772, -9.7422]\times10^3$ |
| $T_F$ | KPCA | 0.2094 | [0.2093, 0.2096] |
|  | KAE | 0.2105 | [0.2103, 0.2106] |
|  | KFDR-B | 0.2100 | [0.2098, 0.2101] |

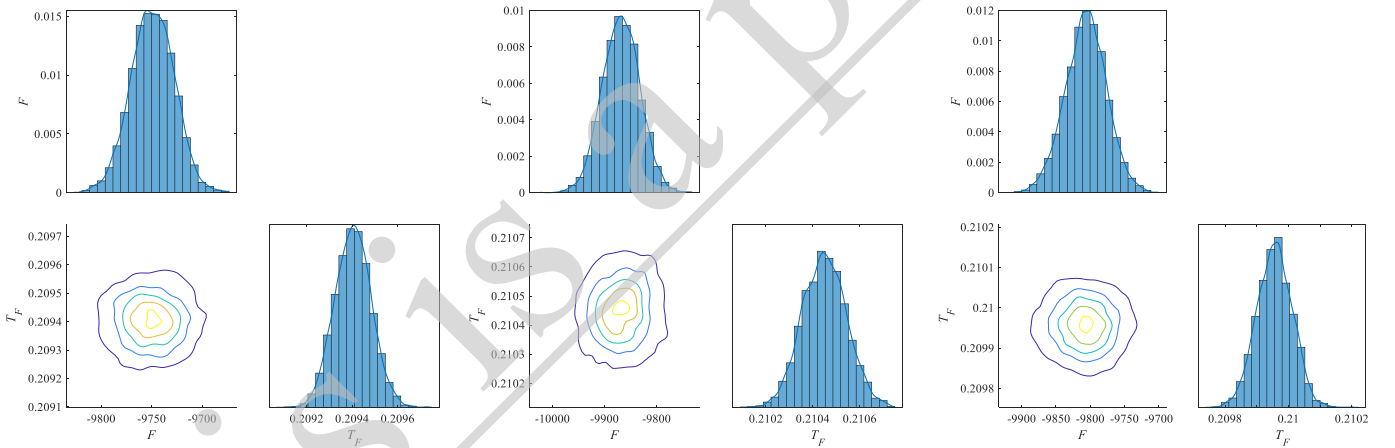

**Fig. 19**. Posterior distributions of the two calibration parameters for the crane structure problem: KPCA model (left), KAE model (middle), and proposed KFDR-B model (right).

## 5. Conclusions and outlook

In this research, we propose a method, referred to as KFDR, that integrates dimension reduction and Kriging surrogate modeling in functional space to perform forward and inverse uncertainty quantification accurately and efficiently for dynamical systems. The proposed KFDR begins by projecting the responses of dynamical systems onto a functional space spanned by a set of predefined basis functions. Next, the functional eigenequation is solved to identify key latent functions, mapping the response of the dynamical system into a low-dimensional latent functional space. Subsequently, Kriging surrogate models with noise terms are constructed in the latent space, enabling accurate and efficient predictions of dynamical systems. Finally, the surrogate model derived from KFDR is directly employed for efficient forward and Bayesian inverse UQ of the dynamical system. Three numerical examples were investigated, leading to the following conclusions:

- By treating the responses of dynamical systems from a functional perspective, they can be represented as linear combinations of a few key latent functions. This functional approach effectively handles noisy data and captures the nonlinear characteristics of the responses. Additionally, an inverse mapping can be directly established from the latent space to the original output space, enabling efficient predictions.
- Kriging surrogate models with noise terms are constructed in the latent functional space to account for errors arising from limited data and feature mapping. Additionally, the probabilistic predictions provided by Kriging models enable the estimation of prediction uncertainty in the time-variant response, allowing metamodeling uncertainty to be considered during uncertainty quantification.
- The illustrative examples demonstrate that the proposed KFDR approach achieves significantly smaller errors in surrogate modeling. Additionally, the forward UQ and inverse UQ results obtained using KFDR show closer agreement with those of the real model compared to the results from the comparative methods, highlighting the accuracy of the proposed approach. Furthermore, the results indicate that B-spline basis functions exhibit greater applicability than Fourier basis functions, making them the recommended choice for the KFDR method.

In the current framework, the Kriging technique is utilized to train surrogate models in the latent functional space. As a result, the proposed method may not be well-suited for high-dimensional inputs, limiting its applicability in scenarios involving complex input excitation functions. These functions typically require a large number of random variables for representation after spectral decomposition. A potential direction for future work is to combine the proposed method with input dimension reduction techniques to broaden its practical applications. Furthermore, the proposed KFDR method is not restricted to uncertainty quantification but can also be extended to reliability analysis and design optimization for dynamical systems. Additionally, with the prediction uncertainty provided by KFDR, it could be integrated with active learning techniques to sequentially enhance the accuracy of the surrogate model.

## CRediT authorship contribution statement

**Zhouzhou Song:** Conceptualization, Methodology, Software, Validation, Writing - original draft, Writing - review & editing, Funding acquisition. **Weiyun Xu:** Software, Writing - review & editing. **Marcos A. Valdebenito:** Supervision, Writing - review & editing. **Matthias G.R. Faes:** Supervision, Funding acquisition, Writing - review & editing.

## Declaration of competing interest

The authors declare that they have no known competing financial interests or personal relationships that could have appeared to influence the work reported in this paper.

## Acknowledgments

The study is partially supported by the Alexander von Humboldt Foundation for the postdoctoral grant of Zhouzhou Song, and the Henriette Herz Scouting program (Matthias G.R. Faes). The authors gratefully acknowledge the supports.

## References

[1] S. L. Brunton, J. L. Proctor and J. N. Kutz, Discovering governing equations from data by sparse identification of nonlinear dynamical systems, Proceedings of the National Academy of Sciences 113 (2016) 3932-3937.

[2] B. Xu, J. He and S. F. Masri, Data-Based Model-Free Hysteretic Restoring Force and Mass Identification for Dynamic Systems, Computer-Aided Civil and Infrastructure Engineering 30 (2015) 2-18.

[3] P. P S, S. Bhartiya and R. D. Gudi, Modeling and Predictive Control of an Integrated Reformer–Membrane–Fuel Cell–Battery Hybrid Dynamic System, Industrial & Engineering Chemistry Research 58 (2019) 11392-11406.

[4] R. Bostanabad, B. Liang, J. Gao, W. K. Liu, J. Cao, D. Zeng, X. Su, H. Xu, Y. Li and W. Chen, Uncertainty quantification in multiscale simulation of woven fiber composites, Computer Methods in Applied Mechanics and Engineering 338 (2018) 506-532.

[5] M. G. R. Faes, M. Daub, S. Marelli, E. Patelli and M. Beer, Engineering analysis with probability boxes: A review on computational methods, Structural Safety 93 (2021) 102092.

[6] D. Zhang, P. Zhou, C. Jiang, M. Yang, X. Han and Q. Li, A stochastic process discretization method combing active learning Kriging model for efficient time-variant reliability analysis, Computer Methods in Applied Mechanics and Engineering 384 (2021) 113990.

[7] Z. Song, H. Zhang, Q. Zhai, B. Zhang, Z. Liu and P. Zhu, A dimension reduction-based Kriging modeling method for high-dimensional time-variant uncertainty propagation and global sensitivity analysis, Mechanical Systems and Signal Processing 219 (2024) 111607.

[8] J. L. Beck and L. S. Katafygiotis, Updating models and their uncertainties. I: Bayesian statistical framework, Journal of Engineering Mechanics 124 (1998) 455-461.

[9] L. S. Katafygiotis and J. L. Beck, Updating models and their uncertainties. II: Model identifiability, Journal of Engineering Mechanics 124 (1998) 463-467.

[10] G. E. P. Box, G. M. Jenkins, G. C. Reinsel and G. M. Ljung, Time Series Analysis: Forecasting and Control, Wiley, 2015.

[11] S. A. Billings, Nonlinear system identification: NARMAX methods in the time, frequency, and spatio-temporal domains, John Wiley & Sons, Chichester, West Sussex, United Kingdom, 2013.

[12] A. J. Smola and B. Schölkopf, A tutorial on support vector regression, Statistics and Computing 14 (2004) 199-222.

[13] V. Ranković, N. Grujović, D. Divac and N. Milivojević, Development of support vector regression identification model for prediction of dam structural behaviour, Structural Safety 48 (2014) 33-39.

[14] D. Xiu and G. E. Karniadakis, The Wiener--Askey Polynomial Chaos for Stochastic Differential Equations, SIAM Journal on Scientific Computing 24 (2002) 619-644.

[15] C. V. Mai, M. D. Spiridonakos, E. N. Chatzi and B. Sudret, Surrogate modeling for stochastic

dynamical systems by combining nonlinear autoregressive with exogenous input models and polynomial chaos expansions, International Journal for Uncertainty Quantification 6 (2016)

[16] C. E. Rasmussen and C. K. I. Williams, Gaussian Processes for Machine Learning, The MIT Press, Cambridge, MA, 2005.

[17] K. Worden, W. E. Becker, T. J. Rogers and E. J. Cross, On the confidence bounds of Gaussian process NARX models and their higher-order frequency response functions, Mechanical Systems and Signal Processing 104 (2018) 188-223.

[18] C. M. Bishop, Neural Networks for Pattern Recognition, Oxford University Press, New York, 1995.

[19] R. F. Mustapa, N. Y. Dahlan, A. I. M. Yassin and A. H. M. Nordin, Quantification of energy savings from an awareness program using NARX-ANN in an educational building, Energy and Buildings 215 (2020) 109899.

[20] S. Schär, S. Marelli and B. Sudret, Feature-centric nonlinear autoregressive models, arXiv preprint (2024) arXiv:2410.07293.

[21] K. Cheng, I. Papaioannou, M. Lyu and D. Straub, State Space Kriging model for emulating complex nonlinear dynamical systems under stochastic excitation, (2024) arXiv:2409.02462.

[22] S. Schär, S. Marelli and B. Sudret, Emulating the dynamics of complex systems using autoregressive models on manifolds (mNARX), Mechanical Systems and Signal Processing 208 (2024) 110956.

[23] E. Jacquelin, N. Baldanzini, B. Bhattacharyya, D. Brizard and M. Pierini, Random dynamical system in time domain: A POD-PC model, Mechanical Systems and Signal Processing 133 (2019) 106251.

[24] G. Blatman and B. Sudret, Sparse polynomial chaos expansions of vector-valued response quantities, CRC Press/Balkema, 2013.

[25] M. Guo and J. S. Hesthaven, Data-driven reduced order modeling for time-dependent problems, Computer Methods in Applied Mechanics and Engineering 345 (2019) 75-99.

[26] M. Vohra, P. Nath, S. Mahadevan and Y.-T. Tina Lee, Fast surrogate modeling using dimensionality reduction in model inputs and field output: Application to additive manufacturing, Reliability Engineering & System Safety 201 (2020) 106986.

[27] Y. Ji, H. Liu, N.-C. Xiao and H. Zhan, An efficient method for time-dependent reliability problems with high-dimensional outputs based on adaptive dimension reduction strategy and surrogate model, Engineering Structures 276 (2023) 115393.

[28] K. Lee and K. T. Carlberg, Model reduction of dynamical systems on nonlinear manifolds using deep convolutional autoencoders, Journal of Computational Physics 404 (2020) 108973.

[29] T. Simpson, N. Dervilis and E. Chatzi, Machine Learning Approach to Model Order Reduction of Nonlinear Systems via Autoencoder and LSTM Networks, Journal of Engineering Mechanics 147 (2021) 04021061.

[30] K. Kontolati, D. Loukrezis, D. G. Giovanis, L. Vandanapu and M. D. Shields, A survey of unsupervised learning methods for high-dimensional uncertainty quantification in black-box-type problems, Journal of Computational Physics 464 (2022) 111313.

[31] A. Lye, A. Cicirello and E. Patelli, Sampling methods for solving Bayesian model updating problems: A tutorial, Mechanical Systems and Signal Processing 159 (2021) 107760.

[32] J. Mercer, Functions of positive and negative type and their connection with the theory of integral equations, Philosophical Transactions of the Royal Society of London Series A 209 (1909) 415-446.

[33] F. O'Sullivan, A Statistical Perspective on Ill-Posed Inverse Problems, Statistical Science 1 (1986) 502-518.

[34] P. Craven and G. Wahba, Smoothing noisy data with spline functions, Numerische Mathematik 31 (1978) 377-403.

[35] C. De Boor, A Practical Guide to Splines, Springer New York, 1978.

[36] Z. Wan, J. Chen, W. Tao, P. Wei, M. Beer and Z. Jiang, A feature mapping strategy of metamodelling for nonlinear stochastic dynamical systems with low to high-dimensional input uncertainties, Mechanical Systems and Signal Processing 184 (2023) 109656.

[37] Y. Liu, L. Li and S. Zhao, A global surrogate model for high-dimensional structural systems based on partial least squares and Kriging, Mechanical Systems and Signal Processing 164 (2022) 108246.

[38] Y. Guo, S. Mahadevan, S. Matsumoto, S. Taba and D. Watanabe, Investigation of Surrogate Modeling Options with High-Dimensional Input and Output, AIAA Journal 61 (2023) 1334-1348.

[39] R. Tripathy, I. Bilionis and M. Gonzalez, Gaussian processes with built-in dimensionality reduction: Applications to high-dimensional uncertainty propagation, Journal of Computational Physics 321 (2016) 191-223.

[40] W. Luo and B. Li, Combining eigenvalues and variation of eigenvectors for order determination, Biometrika 103 (2016) 875-887.

[41] Z. Song, Z. Liu, H. Zhang and P. Zhu, An improved sufficient dimension reduction-based Kriging modeling method for high-dimensional evaluation-expensive problems, Computer Methods in Applied Mechanics and Engineering 418 (2024) 116544.

[42] J. Goodman and J. Weare, Ensemble samplers with affine invariance, Communications in applied mathematics and computational science 5 (2010) 65-80.

[43] S. Marelli and B. Sudret, UQLab: A Framework for Uncertainty Quantification in Matlab, Proc. 2nd Int. Conf. on Vulnerability, Risk Analysis and Management (ICVRAM2014), Liverpool, United Kingdom (2014) 2554-2563.

[44] P.-R. Wagner, J. Nagel, S. Marelli and B. Sudret, UQLab user manual - Bayesian inference for model calibration andinverse problems, Report UQLab-V2.1-113, Chair of Risk, Safety and Uncertainty Quantification, ETH Zurich, Switzerland, 2024.